\documentclass[11pt]{amsart}

\usepackage{hyperref,amsmath,amsthm,amsfonts,amscd,flafter,epsf,amssymb}
\usepackage{epsfig}
\usepackage{amsfonts}
\usepackage{amssymb}
\usepackage{amsmath}
\usepackage{amsthm}
\usepackage{latexsym}
\usepackage{amscd}
\usepackage{flafter}
\usepackage{epsf}
\newcommand{\Q}{\mathbb{Q}}

\newcommand{\xfrak}{\mathfrak{x}}
\newcommand{\afrak}{\mathfrak{a}}
\newcommand{\bfrak}{\mathfrak{b}}
\newcommand{\ifrak}{\mathfrak{i}}
\newcommand{\Ll}{\mathbf{L}}
\newcommand{\tpi}{\widetilde{\pi}}
\newcommand{\tSig}{\widetilde{\Sig}}
\newcommand{\Real}{\mathrm{Re}}
\newcommand{\trace}{\mathrm{tr}}
\newcommand{\Index}{\mathrm{Index}}

\newcommand{\Fbar}{\overline{F}}

\newcommand{\MapModbar}{\overline{\mathcal{X}}}

\newcommand{\Sig}{\Sigma}
\newcommand{\ra}{\rightarrow}
\newcommand{\xra}{\xrightarrow}
\newcommand{\Ht}{\mathrm{H}}

\newcommand{\m}{\mathfrak{m}}
\newcommand{\sig}{\sigma}

\newcommand{\Hbb}{\mathbb{H}}

\newcommand{\C}{\mathbb{C}}

\newcommand{\Qbb}{\mathbb{Q}}

\newcommand{\Rbb}{\mathbb{R}}
\newcommand{\RR}{\Rbb}
\newcommand{\R}{\Rbb}

\newcommand{\Z}{\mathbb{Z}}

\newcommand{\Bcal}{\mathcal{B}}

\newcommand{\Hcal}{\mathcal{H}}
\newcommand{\Ical}{\mathcal{I}}
\newcommand{\Jcal}{\mathcal{J}}
\newcommand{\Kcal}{\mathcal{K}}
\newcommand{\Lcal}{\mathcal{L}}
\newcommand{\Mcal}{\mathcal{M}}

\newcommand{\Ocal}{\mathcal{O}}

\newcommand{\Tcal}{\mathcal{T}}

\newcommand{\Ycal}{\mathcal{Y}}

\newcommand{\ov}{\overline}

\newcommand\CP[1]{{\mathbb{CP}}^{#1}}

\newcommand{\dbar}{\ov{\partial}}

\def\endproof{\relax\ifmmode\expandafter\endproofmath\else
  \unskip\nobreak\hfil\penalty50\hskip.75em\hbox{}\nobreak\hfil\bull
  {\parfillskip=0pt \finalhyphendemerits=0 \bigbreak}\fi}
\def\endproofmath$${\eqno\bull$$\bigbreak}
\def\bull{\vbox{\hrule\hbox{\vrule\kern3pt\vbox{\kern6pt}\kern3pt\vrule}\hrule}}

\newcommand{\Ker}{\mathrm{Ker}}

\newcommand{\Coker}{\mathrm{Coker}}

\newcommand{\Image}{\mathrm{Im}}

\newcommand{\Ind}{\mathrm{Index}}
\newcommand{\Hom}{\mathrm{Hom}}

\newtheorem{thm}{Theorem}[section]
\newtheorem{prop}[thm]{Proposition}

\newtheorem{lem}[thm]{Lemma}

\newtheorem{defn}[thm]{Definition}

\newtheorem{remark}[thm]{Remark}

\newcommand{\Mod}{\Mcal}
\newcommand{\Modbar}{\ov{\Mcal}}

\newcommand{\Fbundle}{\mathcal{F}}

\newcommand{\Group}{\mathfrak{G}}

\newcommand{\MapMod}{\mathcal{X}}

\newcommand{\lra}{\longrightarrow}

\newcommand{\Bundle}{\mathcal{E}}

\newcommand{\Nod}{\mathcal{N}}

\newcommand{\PP}{\mathbb{P}}
\newcommand{\Fbundlebar}{\overline{\mathcal{F}}}

\begin{document}
\title[On finiteness and rigidity in symplectic $3$-folds]
{On finiteness and rigidity of $J$-holomorphic curves in symplectic three-folds }%
\author{Eaman Eftekhary}%
\address{School of Mathematics, Institute for Research in Fundamental Science (IPM),
P. O. Box 19395-5746, Tehran, Iran}%
\email{eaman@ipm.ir}
\begin{abstract}
Given a symplectic three-fold $(M,\omega)$  we
show that for a generic almost complex structure $J$ 
which is compatible with $\omega$,
there are finitely many $J$-holomorphic curves in $M$ of any genus $g\geq 0$
representing a homology class $\beta$ in $\Ht_2(M,\Z)$
with $c_1(M).\beta=0$, provided that the divisibility of $\beta$ is at most $4$ (i.e.
if $\beta=n\alpha$ with $\alpha\in\Ht_2(M,\Z)$ and $n\in \Z$ then $n\leq 4$).
Moreover, each such curve is embedded and $4$-rigid.
\end{abstract}
\maketitle
\section{Introduction}
Let $(M,\omega)$ be a symplectic three-fold, and $J$ be an element in the space 
$\Jcal^\infty(M,\omega)$ of
smooth almost complex structure on $M$ which are compatible with $\omega$. 
For any given homology class $\beta\in \Ht_2(M,\Z)$ and any given 
genus $g\geq 0$, the virtual dimension of 
the moduli space of $J$-holomorphic curves representing $\beta$ which 
have genus $g$ is equal to $c_1(M).\beta$.
In particular, when $c_1(M).\beta=0$ (e.g. if $c_1(M)=0$) this 
moduli space is expected to be zero-dimensional. 
In fact, a conjecture of Ionel and Parker predicts that for a 
generic almost complex structure $J$ all such 
moduli spaces are compact zero dimensional manifolds 
(see section 7.4 of \cite{CK}). In this paper we will 
present a proof when the homology class $\beta$ is not divisible by large integers.

\begin{defn}
 For a homology class $\beta\in\Ht_2(M,\Z)$ define the
  {\emph{divisibility}} $|\beta|$ 
to be the largest integer $n$ such that $\beta=n\alpha$ for some 
$0\neq \alpha\in\Ht_2(M,\Z)$.
\end{defn}

The following is the main result of this paper:
\begin{thm}\label{thm:main-introduction}
For any $J$ in  a subset 
$\Jcal_{rigid}^\infty(M,\omega)\subset \Jcal^\infty(M,\omega)$
of second category, any given genus $g\geq 0$, and any homology class
$\beta\in\Ht_2(M,\Z)$ which satisfies 
$c_1(M).\beta=0$ and $|\beta|\leq 4$, the moduli space
 $\Mod_g(M,\beta;J)$ of somewhere injective  $J$-holomorphic curves of genus $g$
 representing the homology class $\beta$ consists of finitely many elements. 
 Moreover, every curve in this moduli space is embedded and $4$-rigid in the 
 symplectic category.
\end{thm}
An embedded $J$-holomorphic curve $C\subset M$ with holomorphic
normal bundle $N_C$ is called $n$-\emph{rigid} (in symplectic category)
if for any holomorphic branched covering map $\pi:\Sig\ra C$ 
from a smooth Riemann surface $\Sig$ to $C$,
of degree less than or equal to $n$, there is no no-zero section 
$X\in \Gamma(\Sig,\pi^*N_C)$ satisfying the 
linearized Cauchy-Riemann equation
\begin{equation}\label{eq:CR-equation}
\nabla X+J\nabla_{j_\Sig}X+(\nabla_X J)d\pi j_{\Sig}=0,
\end{equation}
where $\nabla$ is the covariant derivative associated 
with a metric connection, $j_\Sig$ is the complex structure on $\Sig$,
and the equation takes place in the vector space 
$\Gamma(\Sig,\Omega^{0,1}_\Sig\otimes_J\pi^*N_C)$.
In particular, such a curve is called {\emph{rigid}} if it is $1$-rigid, and is called 
{\emph{super-rigid}} if it is $n$-rigid for all $n\in\Z^{+}$.
If the Hermitian connection corresponding to $\nabla$ is so that the parallel 
translation commutes with $J$, equation~\ref{eq:CR-equation} 
may be re-written in terms of the 
Nijenhuis tensor $N_J$ for $J$ as 
\begin{equation}
 \label{eq:CR-Nijenhuis}
\nabla X+J\nabla_{j_\Sig}X+\frac{1}{4}N_J\left(\partial_J(\imath_C),X\right)=0,
\end{equation}
where $\imath_C$ denotes the embedding of $C$ in $M$.   
Thus for integrable $J$, the notions of super-rigidity from \cite{R-Jim} and symplectic 
super-rigidity agree.\\

The $n$-rigidity of a $J$-holomorphic curve $C\subset M$ implies that for any genus 
$g\geq 0$ and any homology class $\beta\in\Ht_2(M,\Z)$ which is of the form $d[C]$
for some integer $0<d\leq n$, the compactified moduli space $\Modbar_g(M,\beta;J)$ 
(of  $J$-holomorphic curves $f:\Sig\ra M$ of genus $g$ 
representing the homology class $\beta$)
has an open component which may be identified with the moduli space 
$\Modbar_g(C,d[C])$ of genus $g$, degree $d$ branched covers of $C$.
Thus the contribution of $C$ to the Gromov-Witten invariants 
$N_g(M,\beta)$ is well-defined. This contribution will be denoted by $C_g(C,d;J)$. 
For instance, $C_0(C,d;J)$ is equal to $1/d^3$ when $J$ is an 
integrable complex structure and 
$C=\CP1$ is a rational curve with normal bundle $\Ocal(-1)\oplus \Ocal(-1)$.
When $J$-holomorphic curves are super-rigid, global
 Gromov-Witten theory of $M$ is thus reduced 
to the local Gromov-Witten theory of such $J$-holomorphic curves.\\

The significance of super-rigid curves in the holomorphic 
category was observed by Bryan and
Pandharipande \cite{R-Jim} in the study of a local version 
of Gopakumar-Vafa conjecture 
 which describes the Gromov-Witten invariants of a Calabi-Yau three-fold $M$
in terms of (not mathematically defined) integer valued
invariants, called the \emph{Gopakumar-Vafa} invariants \cite{Go-Va-1}.
It was observed in \cite{P} by Pandharipande that the 
contribution of a super-rigid holomorphic curve
to the Gromov-Witten invariants of higher genus 
(corresponding to the multiples of the homology class the
 curve represents) is
independent of its normal bundle. Thus, for an \emph{ideal} 
Calabi-Yau three-fold (where for an integrable
almost complex structure $J$ and associated with each
homology class and genus there are only finitely many 
$J$-holomorphic curves and all of them are embedded
and super-rigid) the Gopakumar-Vafa conjecture is equivalent 
to its local version.\\

Using the spectral flow from the Dolbeault $\dbar$ operator to the 
linearized Cauchy-Riemann operator 
defined by the left hand side of equation~\ref{eq:CR-equation}, 
one may assign a sign in $\{+1,-1\}$ to each element 
of $\Mod_g(\beta;J)$ for $\beta\in\Ht_2(M,\Z)$ with $c_1(M).\beta=0$, provided that 
$J\in\Jcal^\infty_{rigid}(M,\omega)$ and $|\beta|\leq 4$. Counting 
such curves with the assigned sign we obtain the integer numbers $e_g(\beta;J)$.
It is an interesting question to investigate the dependence of $e_g(\beta;J)$ on the
almost complex structure $J$, and possible evaluation of the local 
contributions $C_{g}(C,d;J)$, at least when $d<5$. 
Note that with integrability assumption on $J$, 
the local contributions $C_{g}(C,d;J)$ are independent of 
the normal bundle, and may be denoted by 
$C_{g}(C,d)$ (see \cite{R-Jim}, also section~\ref{sec:conclusion}).
If a similar conclusion could be extended to arbitrary 
$J\in\Jcal^\infty_{rigid}(M,\omega)$,
the independence of $e_g(\beta;J)$ from $J$ would have been an immediate corollary,
using a M\"obius inversion formula as in \cite{R-Jim}.
However, the work of Taubes in dimension two \cite{Taubes} 
suggests that a wall-crossing phenomena 
may appear in this case, forcing non-trivial dependence of $e_g(\beta,J)$ 
on $J\in \Jcal_{rigid}^\infty(M,\omega)$.\\

Another obvious open direction for further investigation is 
the case of curves with higher divisibility of 
the associated homology class $\beta\in\Ht_2(M,\Z)$, 
and the question of finiteness and (super)-rigidity
for such objects.  The main technical issue in extending 
the results of the current paper is the following. 
If a smooth $J$-holomorphic curve $C\subset M$ 
is not super-rigid, there is a Riemann surface $\Sig$ 
admitting an action of a finite group $\Group$ which gives a branched covering map
$\pi:\Sig\ra \Sig/\Group=C$ with the property that  
the pull-back of the normal bundle $N_C$ of 
$C$ to $\Sig$ via $\pi$ admits a non-trivial section satisfying 
equation~\ref{eq:CR-equation}.
In general, $\Group$ may admit large irreducible real 
representations that damage the transversality argument
of section~\ref{sec:transversality}, which is the heart of our proof. 
Although this obstacle seems hard to overcome, 
the author hopes that the techniques and the setup 
used in this paper may be applied in
other moduli problems.\\  

{\bf{Acknowledgements.}}  The author would like to thank 
Gang Tian for his continuous  advice and support,
and for sharing his insight.

\section{Moduli space of embedded Riemann surfaces}\label{sec:moduli-space}
Fix the integers $k,p>1$ and $\ell\geq k$.
We will sometimes drop these integers from our notation for simplicity, and will 
mention them only in the final statements.\\
 
Let $(M,\omega)$ be a compact symplectic three-fold and
fix a second homology class $\beta\in\Ht_2(M,\Z)$.
Let $\Modbar_g$ denote
 the  Deleign-Mumford compactification of the moduli space $\Mod_g$
 of  Riemann surfaces of genus $g$. Let
$\MapModbar$ be the fiber bundle over $\Modbar_{g}$
 whose fiber over $[\Sig]$ is given by
$$\MapModbar_{\Sig}:=\{f:\Sig\ra M\ |\ 
f\in W^{k,p}(\Sig,M)\ \&\ f_*[\Sig]=\beta\},$$
and let $\MapMod$ be the union of fibers over the moduli space $\Mod_g$
of smooth curves.
Here $W^{k,p}(\Sig,M)$ denotes the Sobolev manifold of maps of type $W^{k,p}$ 
from the surface $\Sig$ to the manifold $M$.
Denote the subspace of $\MapMod$ consisting of the somewhere injective
maps by $\MapMod^\circ$. Let $\Jcal=\Jcal^\ell(M,\omega)$ 
be the $C^\ell$ completion of \
the space of  smooth almost complex
structures on $M$ which are compatible with $\omega$.
Let $\Bundle$ be the vector bundle over $\Ycal=\MapMod\times \Jcal$ whose
fiber over the pair $(f:\Sig\ra M,J)$ is the vector space
$$\Bundle_f:=\Gamma^{k-1,p}(\Sig,\Omega^{0,1}_{j_\Sig}\otimes_{J} f^*T'M),$$
where $j_\Sig$ is the complex structure on $\Sig$.
Here $\Omega^{0,1}_{j_\Sig}$ denotes the vector space of $(0,1)$-forms
on   $\Sig$ determined by the complex structure $j_\Sig$, and 
$\Gamma^{k-1,p}$ denotes the space of sections of type $W^{k-1,p}$.
The almost complex structure $J$
defines a section $\dbar_{J}$ of the vector bundle 
$\Bundle|_{\MapMod\times\{J\}}$ over
$\MapMod\times\{J\}$, given by
$\dbar_{J}(f)=df+J\circ df \circ j_{\Sig}$.
The sections $\dbar_{J}$  glue together and define a smooth section 
$\dbar$ of $\Bundle\ra \Ycal$
given by $\dbar(f,J)=\dbar_{J}(f)$.  The following lemma is proved in \cite{OZ}
(see the proof of theorem 1.4):
\begin{lem}\label{lem:1}
With the above notation fixed, the intersection of $\dbar$ with the zero section of
$\Bundle\ra \Ycal$ is transverse.
\end{lem}
We define:
\begin{displaymath}\begin{split}
\Mod_g(\beta,J)&:=\Big\{f:\Sig\ra M\ \Big|\
{\dbar_{J}(f)=0,\
f\in \MapMod^\circ}
\Big\}\subset \MapMod^\circ\times \{J\},\ \ \&\\
\Mod_g(\beta)&:=\Big\{(f,J)\in \Ycal\ \Big|\ \dbar_J(f)=0\Big\}
=\bigcup_{J\in\Jcal}\Mod_g(\beta;J)\subset \Ycal.
\end{split}
\end{displaymath}
Thus the moduli space $\Mod_g(\beta)$
is a smooth, Banach manifold.  The
projection map $\Ycal=\MapMod\times \Jcal\ra \Jcal$ induces a map
$\pi_{g,\beta}:\Mod_g(\beta)\ra \Jcal$, which will be 
Fredholm of index $2c_1(M).\beta$.
By Sard-Smale theorem, the set of regular values for the projection
map $\pi_{g,\beta}:\Mod_g(\beta)\ra \Jcal$ is a subset 
$\Jcal_{reg}(g,\beta)\subset \Jcal$, which is
of second category. For every $J\in \Jcal_{reg}(g,\beta)$,
$\Mod_g(\beta,J)=\pi_{g,\beta}^{-1}(J)$ 
 is a smooth manifold whose dimension
is equal to the index of the projection map $\pi_{g,\beta}$.  By the regularity
of   $J$,  $\dbar_{J}:\MapMod^\circ\times\{J\}\ra \Bundle$
intersects the zero section transversely over $\Mod_g(\beta,J)$. Set
$$\Jcal_{reg,0}:=\bigcap_{g\geq 0}
\bigcap_{\substack{\beta\in\Ht_2(M,\Z)\\ c_1(M).\beta=0}}
\Jcal_{reg}(g,\beta).$$
Let $\Tcal_J$ be the tangent space to $\Jcal$ at $J$,
consisting of the linear homomorphisms $u:TM\ra TM$ which satisfy
$$u\circ J+J\circ u=0\ \ \ \&\ \ \ 
\omega_x(u(X),Y)+\omega_x(Y,u(Y))=0,\ \ \forall X,Y\in T_xM.$$
 The derivative of $\dbar_J$ at a zero $(f,J)$ of this section is a linear map
$$d\dbar:T_f\MapMod\oplus \Tcal_J\ra T_{((f,J);0)}\Bundle=
T_f\MapMod\oplus \Tcal_J\oplus \Bundle_{f}.$$
Projection over the last factor in this decomposition,
composed with the differential $d\dbar$, gives a 
linear operator (after multiplication by 2)
$$\Ll:T_f\MapMod\oplus \Tcal_J=
\Gamma^{k,p}(f^*TM)\oplus \Ht^1(T_\Sig)\oplus \Tcal_J\ra
 \Bundle_f=\Gamma^{k-1,p}(\Omega^{0,1}_\Sig\otimes_J f^*TM).$$
The intersection of $\dbar$ and the zero section at $(f,J)$ is transverse
if and only if this linear map $\Ll$ is surjective.
If $\nabla$ denotes the Levi-Civita connection of the metric on $M$
(defined by $\langle X,Y\rangle=w(X,JY)$), we may write down an explicit 
formula for this linear map
\begin{equation}\label{eq:L&K}
\begin{split}
\Ll(X,\eta,u)&=\nabla X+J \nabla_j X+(\nabla_X J) df  j
+J df \eta+u df j\\
&=:\Lcal(X,\eta)+u df j =:\Lcal(X,\eta)+L(u).\\
\end{split}
\end{equation}
Note that $\Lcal=\Lcal_{(f,J)}$ is in fact a linear operator
$$\Lcal:\Gamma^{k,p}(f^*TM)\oplus \Ht^1(T_\Sig)\lra \Gamma^{k-1,p}
(\Omega^{0,1}_{\Sig}\otimes_J f^*T'M),$$
which is Fredholm of index $2c_1(M).\beta$. By integrating the
point-wise Hermitian inner product of $\Omega^{0,1}_\Sig
\otimes f^*T'M$ we obtain a Hermitian inner product on
 $\Gamma^{k-1,p}(\Omega^{0,1}_\Sig\otimes f^*T'M)$. 
 Similarly, we have an inner product
on $\Gamma^{k,p}(f^*TM)\oplus \Ht^1(T_\Sig)$. Using these inner products
 we may define an adjoint operator for $\Lcal$
$$\Lcal^*:\Gamma^{k,q}(\Omega^{0,1}_\Sig\otimes_J f^*T'M) \lra
\Gamma^{k-1,q}(f^*TM)\oplus \Ht^1(T_\Sig),$$
where $q$ satisfies $\frac{1}{p}+\frac{1}{q}=1$. It follows from 
elliptic regularity that every section
$\delta\in\Gamma^{k-1,p}(\Sig,\Omega_\Sig^{0,1}\otimes_J f^*T'M)$ 
which is annihilated by the image of the operator $\Lcal$ in the sense that
$$\langle \Lcal(X,\eta),\delta\rangle=0,\ \ \ \forall \ X\in
\Gamma^{k,p}(\Sig, f^*T'M),\ \eta\in\Ht^1(\Sig,T_\Sig),$$
is automatically in the kernel of $\Lcal^*$, and is of class 
$W^{k,q'}$ for all values of $q'$. Thus, the cokernel of $\Lcal$ may be 
identified with the kernel of $\Lcal^*$.
Also, note that if $f:\Sig\ra M$ is in $\Mod_g(\beta,J)$, then it is at 
least of class $C^{\ell+1}$ by 
elliptic regularity, and the moduli space $\Mod_g(\beta,J)$ is thus 
independent  of $k$.\\

The following is one of the main results of \cite{OZ}  (theorem 1.4):
\begin{prop}\label{prop:no-node}
For any $J$ in a subset 
$\Jcal_{reg}^{\ell}(M,\omega)\subset 
\Jcal^{\ell}(M,\omega)$ 
of second category the following is true for any 
$\beta\in\Ht_2(M,\Z)$. If $c_1(M).\beta<0$ the moduli 
space  $\Mod_g(J,\beta)$ is empty, while for $c_1(M).\beta=0$  
$\Mod_g(J,\beta)$ consists of embeddings. Moreover,
if $c_1(M).\beta=c_1(M).\beta'=0$  and if
$f\in \Mod_g(\beta,J)$ and $f'\in \Mod_{g'}(\beta',J)$ 
have distinct images, then their images are disjoint. 
\end{prop}

Fix $J\in \Jcal_{reg}^\ell(M,\omega)\cap \Jcal^\infty(M,\omega)$.
For the moment, we are not worried about the existence of 
such almost complex structures,
and will return to this issue in section~\ref{sec:regularity}. Let $\beta\in\Ht_2(M,\Z)$ 
be a homology class with $c_1(M).\beta=0$.
The existence of  a sequence of embedded $J$-holomorphic curves
$f_i:\Sig_i\ra M$ in $\Mod_g(\beta,J)$
implies that the domains converge to a possibly singular 
curve $\Sig\in \Modbar_{g}$, and that a subsequence of $f_i$
will converge to a $J$-holomorphic map $f:\Sig \ra M$, 
by Gromov compactness theorem.
The homology class $\alpha\in\Ht_2(M,\Z)$  associated with the map $f$ 
 may be different from $\beta$.
The class $\beta-\alpha$ will be represented by $J$-holomorphic
\emph{bubbles} (i.e. $J$-holomorphic maps from $\CP1$ to $M$).
Thinking of these bubbles as the components of 
$\Sig$ and dropping the stability condition on the
domain, we may assume that $f_*[\Sig]=\beta$.
The normalization of the domain $\Sig$
will be a union $\Sig^1\cup\Sig^2\cup...\cup \Sig^k$ 
of smooth components. The map $f$ induces a $J$-holomorphic
map $f^i:\Sig^i\ra M$ on each component $\Sig^i$ which 
decomposes as $f^i=\imath^i\circ \pi^i$, where $\imath^i$ is a somewhere
injective map from a domain $C^i$ of some genus 
$h_i$ to $M$ and $\pi^i:\Sig^i\ra C^i$ is a branched
 covering (which may be $Id_{\Sig^i}$).
When $\Sig^i$ is a collapsed component, 
we will assume that $C^i=\Sig^i,\ \imath^i:C^i\ra M$
is the collapsing (constant) map, and $\pi^i$ is the identity map of $\Sig^i=C^i$.\\
 
If the homology class represented by $\imath^i$ is $\alpha_i$, $\imath^i$ 
would be in $ \Mod_{h_i}(\alpha_i,J)$. Since this
moduli space is non-empty, by 
proposition~\ref{prop:no-node} we would have $c_1(M).\alpha_i\geq 0$. 
Since $\beta=m_1\alpha_1+...+m_k\alpha_k$ 
for  positive integers $m_i\in \Z^+,\ i=1,...,k$ and $c_1(M).\beta=0$,
we conclude that $c_1(M).\alpha_i=0$ for $i=1,...,k$.
Proposition~\ref{prop:no-node} implies that the images of
$\imath^i$ are disjoint or identical. Since $\Sig$ is connected,  there is a
$J$-holomorphic curve $C$ in $M$ such that each 
map $\imath^i$ is in fact the embedding of $C^i\simeq C$
in $M$. Thus the map $f:\Sig\ra M$ may be regarded 
as the composition of a branched covering map 
$\pi:\Sig\ra C$ (which is equal to $\pi^i$ on the 
component $\Sig^i$) with the $J$-holomorphic
embedding $\imath_C:C\ra M$.\\

A neighbourhood of $C$ in $M$ may be (symplectically) 
identified with a neighbourhood of the zero section
in its normal bundle $N_C$. For $m$ sufficiently
large we may assume that the image of $f_m$ is  included in this neighbourhood.
We may thus work in  $N_C$ and assume that the images of all $f_m$ are in the 
unit disk bundle around the zero section in $N_C$ 
(implicitly, we are fixing a metric on the normal
bundle $N_C$).\\

For any constant $\epsilon>0$, let $\chi_\epsilon:\R\ra [0,1]$ 
be a smooth non-decreasing function with $\chi_\epsilon(x)=\epsilon$
for $x\leq 1$ and $\chi_\epsilon(x)=1$ for $x\geq 2$. Furthermore,
 assume that $\chi_\epsilon$ varies smoothly with 
$\epsilon$ for $0\leq \epsilon\leq 1$.
Let $F_\epsilon:N_C\ra N_C$ be the diffeomorphism defined by 
$$F_\epsilon(z;X):=(z;\chi_\epsilon(\|X\|^2).X),\ \ 
\forall\ z\in C\ \&\  X\in (N_C)_z.$$
Let $J_\epsilon$ denote the almost complex structure $F_\epsilon^*J$ on $N_C$.
For any embedding $f_m$ as above, let $\epsilon(m)$ denotes
the supremum norm of $f_m$ as a multi-section of $N_C$. 
The embedding $f_m$ may then be composed with 
$F_{\epsilon(m)}^{-1}$ to obtain a 
$J_{\epsilon(m)}$-holomorphic embedding of $\Sig_m$ in $N_C$, 
which will be denoted by $h_m:\Sig_m\ra N_C$. 
By the above construction, the multi-sections $h_m$ will have 
supremum norm equal to $1$.\\

As the sequence $\epsilon(m)$ converges to zero with $m$ going to infinity, the 
almost complex structures $J_{\epsilon(m)}$ converge to a limit denoted by $J'$. 
Gromov compactness theorem then tells us that the sequence $h_m:\Sig_m\ra N_C$ 
converges to a limit $h:\Sig\ra N_C$, which is $J'$-holomorphic.
The composition of $h$ with the projection of $N_C$ over $C$ is the 
branched covering map $\pi:\Sig\ra C$, and the map $h$ is  determined 
by a section $X\in \Gamma(\Sig,\pi^*N_C)$. 
The section $X$ is in fact a union of sections
$$X^j\in\Gamma^\infty\left(\Sig^j,(\pi^j)^*N_C\right)
\subset\Gamma^{k,p}\left(\Sig^j,(\pi^j)^*N_C\right),\ \ j=1,...,k.$$
The map $h$ is $J'$-holomorphic if and only if for $j=1,...,k$, 
the section $X^j$ satisfies the following equation
$$\Kcal_{\pi^j}(X^j):=\nabla 
X^j+J\nabla_{j_{\Sig^j}}X^j+\left(\nabla_{X^j}J\right)d\pi^jj_{\Sig}=0$$ 
in $\Gamma^{k-1,p}(\Sig^j,\Omega_{\Sig^j}^{0,1}\otimes_J(\pi^j)^*N_C)$.
This means that $X^j$ gives an element in the kernel of the operator
$$\Kcal_{\pi^j}:\Gamma^{k,p}(\Sig^j,(\pi^j)^*N_C)\ra\Gamma^{k-1,p}(\Sig^j,
\Omega^{0,1}_{\Sig^j}\otimes (\pi^j)^*N_C).$$
See \cite{Oh}, theorem 7.1 for a similar proof, also \cite{IP2} proposition 6.6.
The section would be non-zero on at least one of the non-collapsing components 
$\Sig^j$, for $j\in\{1,...,k\}$, since the supremum norm of the 
section $X$ is equal to $1$, and the Riemann surface $\Sig$ is connected.\\

Thus the existence of a sequence as above which converges to an embedded 
$J$-holomorphic curve $C$ of genus $h\leq g$ implies that there is a 
branched covering map $\pi_S:S\ra C$ (with $S=\Sig^j$ and $\pi_S=\pi^j$ for some 
$j$), and a non-trivial section $X$ in the kernel of the elliptic operator
$$\Kcal_{\pi_S}:\Gamma^{k,p}(S,\pi_S^*N_C)\ra
\Gamma^{k-1,p}(S,\Omega_{S}^{0,1}\otimes_J(\pi_S)^*N_C).$$
Based on this observation, we make the following definition, which should 
be compared with the concepts of rigidity and super-rigidity in the 
integrable case by  Bryan and Pandharipande  \cite{R-Jim}.\\

\begin{defn}
\label{defn:rigidity}
Let $J\in\Jcal^\infty(M,\omega)$, and $C$ be a smooth embedded $J$-holomorphic 
curves of genus $h\geq 0$ in $M$ with normal bundle  $N_C$.
The curve $C$ is called $(m,n)$-\emph{rigid} for the integers $m\geq 0$ 
and $n>0$, if for every branched covering map $\pi_S:S\ra C$ from a 
smooth Riemann surface $S$ whose genus $g$ satisfies $g\leq m+h$, 
and with $\mathrm{deg}(\pi_S)\leq n$, the kernel of the linear operator 
\begin{equation}
 \label{eq:operator-K}
\begin{split}
 &\Kcal_{\pi_S}:\Gamma^{k,p}(S,\pi_S^*N_C)\ra
\Gamma^{k-1,p}(S,\Omega_{S}^{0,1}\otimes_J(\pi_S)^*N_C)\\
&\Kcal_{\pi_S}(X):=\nabla X+J\nabla_{j_S}X+\left(\nabla_XJ\right)\circ d\pi_S\circ j_S
\end{split}
\end{equation}
is trivial. For a positive integer $n>0$, the curve $C$ is called $n$-\emph{rigid}
if it is $(m,n)$-rigid for all $m\in\Z^{>0}$. The curve $C$ is called super-rigid if it is 
$n$-rigid for every integer $n>0$.
This definition does not depend on the particular choice of the integers $k,p\geq 2$ 
by elliptic regularity.
In a similar way, we may define the notions of $n$-rigidity and super-rigidity if $J$ 
is an almost complex structure of class $C^\ell$ and the embedding of $C$ in $M$
is of class $W^{k,p}$ for some integers $k,p\geq 2$ such that $k\leq \ell$.   
\end{defn}

Our discussion in this section implies that for an almost complex structure $J$ in 
$\Jcal_{reg}^\ell(M,\omega)\cap\Jcal^\infty(M,\omega)$, 
if an embedded $J$-holomorphic
curve $C\in \Mod_h(M,\alpha,J)$ is the limit of a sequence of embedded
$J$-holomorphic curves $f_i:\Sig_i\ra M$ representing the homology class 
$\beta$ (which is forced to be of the form $n\alpha$) in the sense that these 
curves converge to a multiple cover $\Sig$ of $C$ of some genus $g=h+m$, then 
$C$ is not $(m,n)$-rigid.

\section{Rigidity; a few reductions}\label{sec:reductions}

In this section, we will study a few possible reductions of 
the concept of rigidity, which would be 
useful in our later considerations.\\

{\bf {Reduction 1.}}
Suppose that $\pi:\Sig\ra C$ is an arbitrary branched covering map 
between smooth Riemann surfaces. Suppose that $B\in \mathrm{Div}(\Sig)$
is the branching divisor of $\pi$, and that 
$\pi(B)=\{q_1,...,q_l\}\subset C$ is the set of critical values of $\pi$.
Fix $q\in C-\pi(B)$ and let $p_i,\ i=1,...,n$ denote all the points 
with $\pi(p_i)=q$. The group
$\pi_1(C-\pi(B),q)$ acts on $\pi^{-1}(q)=\{p_1,...,p_n\}$ as follows. Suppose that
$$\gamma:[0,1]\ra C-\pi(B),\ \ \ \gamma(0)=\gamma(1)=q$$
is a loop representing an element $[\gamma]\in\pi_1(C-\pi(B),q)$. 
For $i=1,...,n$, we may find the lift $\gamma_i$ of $\gamma$ under the 
covering map $\pi$ such that
$$\gamma_i:[0,1]\ra \Sig-B,\ \ \gamma_i(0)=p_i,\ \ 
\&\ \ \pi\circ \gamma_i=\gamma.$$
We may then define the action of $[\gamma]$ on $\{p_1,...,p_n\}$ by
$$[\gamma]\star p_i:=\gamma_i(1)\in\pi^{-1}(q)=\{p_1,...,p_n\},\ \ \ \forall\ 
i\in\{1,...,n\}.$$
This gives a homomorphism 
$$\rho:\pi_1(C-\pi(B),q)\lra S_n$$
from the fundamental group of $C-\pi(B)$ to the group 
$S_n$ of permutations in $n$ letters.
The kernel $\Ker(\rho)$ of this homomorphism, 
which is a normal subgroup of the fundamental group, determines a 
finite regular covering map $\pi^\circ:\Sig^\circ \ra C-\pi(B)$, 
which may be extended to a branched covering map
$\tpi:\tSig\ra C$ from a compact Riemann surface $\tSig$ to $C$. 
Moreover, this map decomposes as $\tpi=\pi\circ \tau$ for some
branched covering map $\tau:\tSig\ra \Sig$. Since $\Ker(\rho)$ is 
normal in the fundamental group, 
the group of deck transformations for $\pi^\circ$  may be computed as 
$$\mathrm{Dec}(\pi^\circ)=\frac{\pi_1(C-\pi(B),q)}
{\pi^\circ_*(\pi_1(\Sig^\circ,p^\circ))}
=\frac{\pi_1(C-\pi(B),q)}{\Ker(\rho)}\simeq \Image(\rho)<S_n.$$
Denote this later subgroup of $S_n$ by $\Group$. 
Note that the degree of $\tpi$ is equal to the 
order of the finite group $\Group$. Any deck transformation of 
$\Sig^\circ$ may be extended to $\tSig$ as an
automorphism with possible fixed points, and 
$C$ may be realized as the quotient of $\tSig$ under the 
action of $\Group<\mathrm{Aut}(\tSig)$.\\

Now suppose that the curve $C$ is an embedded 
$J$-holomorphic curve in $M$ for some almost complex 
structure $J\in\Jcal^\infty(M,\omega)$, and let $N_C$ 
denote the normal bundle of $C$. Fix the branched 
covering map $\pi:\Sig\ra C$ and consider the operator
$$\Kcal_\pi:\Gamma^{k,p}(\Sig,\pi^*N_C)\lra \Gamma^{k-1,p}
\left(\Sig,\Omega_\Sig^{0,1}\otimes_J\pi^*N_C\right)$$
introduced earlier. If $X\in\Ker(\Kcal_\pi)$ is  non-trivial,
then $\tau^*X\in\Ker(\Kcal_{\tpi})$ is non-trivial as well. 
Thus, if an embedded $J$-holomorphic curve 
$C$ is not $n$-rigid, there is a Riemann surface 
$S=\tSig$ admitting an action of a subgroup 
$\Group$ of $S_n$, with $\pi_S:S\ra C=S/\Group$ the 
corresponding branched covering map,
so that the kernel of $\Kcal_{\pi_S}$ is non-trivial. \\

{\bf{Reduction 2.}}
Let us now assume that $\pi=\pi_S:S\ra S/\Group=C$ is a branched covering map 
coming from the action of a finite group, and assume that $X\in\Ker(\Kcal_\pi)$
is a non-trivial section. For any $\sig\in\Group$, it is then clear that 
$\sig^*X$ is also in $\Ker(\Kcal_\pi)$. Thus the action of the 
group ring $\R_\Group$ on $\Gamma^{k,p}(S,\pi^*N_C)$ defined by
$$\afrak.Y:=\sum_{sig\in\Group} a_\sig \sig^*Y,\ \ \forall\ 
\afrak=\sum_{\sig\in\Group}a_\sig.\sig^{-1}\in\R_\Group,\  
\&\ Y\in\Gamma^{k,p}(S,\pi^*N_C)$$
induces an action of $\R_\Group$ on $\Ker(\Kcal_\pi)$.
In other words, if $X\in\Ker(\Kcal_\pi)$ and $\afrak\in\R_\Group$, then 
$\afrak.X\in\Ker(\Kcal_\pi)$.
Let us denote by $\m_X$ the left ideal of $\R_\Group$ 
consisting of the elements $\afrak\in\R_\Group$
such that $\afrak.X=0$.
Let $I(\Group)$ be the set of irreducible real representations $\ifrak$ of $\Group$.
The group ring $\R_\Group$ may then be decomposed  
(by Artin-Wedderburn and Maschke theorems) as
\begin{equation}
 \label{eq:decomposition-R_G}
\R_\Group\simeq \bigoplus_{\ifrak\in I(\Group)}M_{\ell(\ifrak)}(R_\ifrak),\ \ 
R_\ifrak\in\{\R,\C,\Hbb\},
\end{equation}
where $M_\ell(R)$ denotes the space of $\ell\times\ell$ matrices with 
entries in the ring $R$, $\ell(\ifrak)$
is the dimension of the representation $\ifrak$ as an algebra over $R_\ifrak$, and 
$\Hbb$ denotes the ring of quaternions. Let us denote by $\imath_\ifrak\in\R_\Group$
the identity matrix in the matrix algebra associated with the representation $\ifrak$.
We thus have
$$1_{\R_\Group}=\sum_{\ifrak\in I(\Group)}\imath_\ifrak\ \ 
\Rightarrow\ \ X=\sum_{\ifrak\in I(\Group)}\imath_\ifrak . X=:
\sum_{\ifrak\in I(\Group)}X_\ifrak.$$
If $X$ is in $\Ker(\Kcal_\pi)$, so are $X_\ifrak$. Moreover,
 at least one of the sections $X_\ifrak$ is non-zero.
Fix such a representation $\ifrak\in I(\Group)$ and let 
$R=R_\ifrak$ and $\ell=\ell(\ifrak)$.
Let $\epsilon_i,\ i=1,...,\ell$ denotes the matrix in $M_\ell(R)$ 
with $1$ as its $(i,i)$ entry and zeros elsewhere.
Note that we may think of $\epsilon_i$ as an element of $\R_\Group$. Since 
$$0\neq X_\ifrak=\imath_\ifrak.X=(\epsilon_1+...+\epsilon_\ell).X,$$
at least one of the sections $\epsilon_i.X$ is a non-zero 
element in $\Ker(\Kcal_\pi)$. Denote one such 
section by $Y$. The left ideal $\m=\m_Y$ is then a maximal
 left ideal $\m_{\ifrak,i}$ which consists of all elements of 
$\R_\Group$ which are characterized using the presentation of 
equation~\ref{eq:decomposition-R_G} as those matrices which
have zeros in the $i$-th column in the matrix presentation
 corresponding to the representation $\ifrak$.
The maximal left ideal $\m$ determines a restriction of the operator $\Kcal_\pi$:
$$\Kcal_\pi^\m:\Gamma_\m^{k,p}(S,\pi^*N_C)\lra \Gamma_\m^{k-1,p}
\left(S,\Omega_S^{0,1}\otimes_J \pi^*N_C\right).$$
Here $\Gamma^\bullet_\m(S,\star)$ denotes the vector space of sections 
$Z$ of the bundle $\star$ such that the left 
ideal $\m_Z$ contains $\m$. The section $Y$ is then in the kernel 
of the operator $\Kcal_\pi^\m$. Note that associated 
with any covering map $\pi:S\ra C=S/\Group$ as above, there are 
finitely many maximal left ideals of the form 
$\m=\m_{\ifrak,i}$. We define $\ell(\m)=\ell(\ifrak)$ and $R_\m=R_\ifrak$.
We  denote the finite set of such maximal left ideals by $I(\pi)$.\\

The kernel $\Ker(\rho)$ of the representation $\rho:\Group\ra M_\ell(R)$ 
is a normal subgroup of $\Group$ and determines a degeneration of the covering 
map $\pi:S\ra C=S/\Group$ as a composition
$$S\xra{\pi_1}S'=S/\Ker(\rho)\xra{\pi_2}C=S/\Group=S'/\Image(\rho).$$
Moreover, since the elements of $\Ker(\rho)$ preserve the section 
$Y$ of the bundle $\pi^*N_C$, $Y$ should be of the form $\pi_1^*Z$ for some 
$Z\in\Gamma^{k,p}(S',\pi_2^*N_C)$. Clearly, $Z\in \Ker(\Kcal_{\pi_2})$ is non-trivial,
since $Y$ is non-trivial. We may thus replace $S$ by $S'$, $\Group$ by
 $\Group'=\Image(\rho)$, and the section $Y$ by $Z$.
The above considerations imply that if an embedded 
$J$-holomorphic curve $C$ is not $n$-rigid, then there is a Riemann surface 
$S$ admitting an action of a finite group $\Group<S_n$ (corresponding to a 
branched covering map $\pi:S\ra C=S/\Group$), and a maximal left ideal 
$\m\in I(\pi)$ of the group ring $\R_\Group$, so that the kernel of the operator
$$\Kcal_\pi^\m:\Gamma_\m^{k,p}(S,\pi^*N_C)\lra 
\Gamma_\m^{k-1,p}\left(S,\Omega_S^{0,1}\otimes_J \pi^*N_C\right)$$
is non-trivial. Moreover, the representation
$$\rho:\Group\ra M_\ell(R), \ \ \ \ell=\ell(\m)\ \& \ R=R_\m\in\{\R,\C,\Hbb\}$$
is faithful.\\

In the above situation, for any $Y\in\Gamma^{k,p}_\m(S,\pi^*N_C)$, 
we introduce the formal sum
$$\alpha(Y):=\sum_{\sig\in\Group}(\sig^*Y).\sig.$$
Any element $\afrak=\sum_{\sig\in\Group}a_\sig.\sig^{-1}$ 
of the group ring $\RR_\Group$ acts 
on any such expression from left and right by
\begin{displaymath}
 \label{eq:action-on-a(Y)}
\begin{split}
 \afrak\odot \alpha(Y):=\sum_{\sig,\tau\in\Group}a_\sig(\tau^*Y).
 \sig^{-1}\tau,\ \ \&\ \ 
\alpha(Y)\odot\afrak:=\sum_{\sig,\tau\in\Group}a_\sig(\tau^*Y).\tau\sig^{-1}.
\end{split}
\end{displaymath}
From the above two definitions we may compute
\begin{displaymath}
 \begin{split}
  &\afrak\odot\alpha(Y)=\sum_{\sig,\tau\in\Group}a_\sig
  \left((\sig\tau)^*Y\right).\tau
=\sum_{\tau\in\Group}\tau^*\left(\afrak.Y\right).Y=\alpha(\afrak.Y),\ \ \&\\
  &\tau^*\alpha(Y)=\tau^*\left(\sum_{\sig\in\Group}\sig^*Y.\sig\right)
=\sum_{\sig\in\Group}\tau^*\sig^*(Y).\sig=\alpha(Y)\odot \tau^{-1}\ \ \ \forall\ 
\tau\in\Group.
\end{split}
\end{displaymath}

{\bf{Reduction 3.}}
In order to study $4$-rigidity, as observed in the first reduction, 
we only need to consider the case where the group $\Group$ is isomorphic to 
a subgroup of $S_4$. The subgroups of $S_4$ are isomorphic to one of
the following groups:
$$\Z/k\Z,\ \ k=1,2,3,4,\ \ \Z/2\Z\oplus\Z/2\Z,\ \ S_3,\ \ D_8,
\ \ A_4\ \ \&\ \ S_4.$$
Except for the last $4$ in this list, the other groups are abelian and the
matrix algebras corresponding to their irreducible real representations 
are either $M_1(\R)$ or $M_1(\C)$.
The group $A_4$ has $3$ real irreducible representations with the 
corresponding matrix
algebras equal to $M_1(\R),M_1(\C)$ and $M_3(\R)$. For $D_8$, the dihedral group, 
we may distinguish $4$ one dimensional irreducible representations with the 
associated matrix algebra equal to $M_1(\R)$, and 
one two dimensional irreducible representation with 
the corresponding matrix algebra equal to $M_2(\R)$. Finally for the groups 
$S_3$ and $S_4$ the associated matrix algebra of any irreducible real representation
is one of $M_1(\R), M_2(\R)$ and $M_3(\R)$. 
These observations imply that if $\Group$ 
is a subgroup of $S_4$ acting on a Riemann surface $S$, 
and if $\pi$ denotes the branched covering map 
$S\ra C=S/\Group$, then for any $\m\in I(\pi)$ the associated 
matrix algebra is one of $M_i(\R),\ i=1,2,3$ or $M_1(\C)=\C$.\\

\section{Moduli space of multiply covered curves}\label{sec:transversality}
Once again, fix the integers $k,p>1$ and $\ell\geq k$ throughout this section.
For $J\in\Jcal_{reg}^\ell(M,\omega)$ let $C\in\Mod_h(M,\beta;J)$ 
be a $J$-holomorphic curve of genus $h$ and class $W^{k,p}$
(and hence automatically of class $W^{\ell+1,p}$) in $M$ which represents a homology 
class $\alpha\in\Ht_2(M,\Z)$ with $c_1(M).\alpha=0$. 
Suppose that $\pi:\Sig\ra C$ is a branched covering map,
and $B\in \mathrm{Div}(\Sig)$ is the branching divisor. 
Let $\{q_1,...,q_l\}$ be the set of points in
the image of $B$ (i.e. the set of critical values for $\pi$). 
When we move these $l$ punctures on $C$, the 
complex structure on the complement of $q_1,...,q_l$ on 
the Riemann surface $C$ determines a corresponding 
complex structure on $\Sig$. This gives a family of branched covering maps 
from Riemann surfaces of genus $g=g_\Sig$ to $C$ with the same monodromy data as 
the original branched covering map $\pi$. Note that some branched covering maps 
may have multiple representatives in this description due to automorphisms.
This determines an open family of covers of $C$, 
which has dimension $l$ (note that the embedded 
curve $C$ is parametrized). Moreover, as $(C,J)$ 
moves in $\Mod_h(\alpha)$ the above open families
glue together and form a moduli space which will 
be denoted by $\Mod_\pi(\beta)$, with 
$\beta=\deg(\pi).\alpha$. Any point in this moduli 
space corresponds to a branched covering map
with the same branching behaviour and monodromy information 
as the initial map $\pi$.
We will say in short that any branched covering 
map in this family has the same {\emph{topological type}}
as $\pi$. Having fixed the target genus $h$, 
the domain genus $g\geq h$ and the degree $n$, 
there are only finitely many topological types of branched covering maps of degree $n$ 
from a smooth domain of genus $g$ to  a smooth 
target of genus $h$.
If the initial branched covering map $\pi$ comes 
from the action of a group $\Group$ on $\Sig$
(so that $\pi:\Sig\ra C=\Sig/\Group$ is the quotient map), 
it is clear that all the maps of the same 
topological type as $\pi$ come from an action of $\Group$, 
since the action may be regarded as the action 
of deck transformations of the honest covering map obtained by 
removing the critical values from
$C$ and their pre-images from $\Sig$. In this situation 
we say that the {\emph{topological type of the group
action}} for other branched covering maps in the family is 
the same as the topological type of the group action 
for $\Group$ and the branched covering map $\pi:\Sig\ra C=\Sig/\Group$.\\

Fix the topological type of the group action and 
the corresponding branched covering map 
$\pi:\Sig\ra C=\Sig/\Group$ as above. Let 
$n=\deg(\pi)=|\Group|$, $\beta=n\alpha\in\Ht_2(M,\Z)$,
and assume that $g$ denotes the genus of $\Sig$.
The moduli space $\Mod_\pi(\beta)$  is then fibered 
over $\Mod_h(\alpha)$, where the fiber over
$(\imath_C:C\ra M,J)$ consists of the $l$ dimensional 
family of branched covering maps with target $C$ which are
of the same topological type as $\pi$. The moduli space 
$\Mod_\pi(\beta)$ is a smooth Banach manifold
and the fibration $q_\pi:\Mod_\pi(\beta)\ra \Mod_h(\alpha)$ 
is a Fredholm map of index $l$.
We will denote the composition of $q_\pi$ with the projection map from 
$\Mod_h(\alpha)$ to $\Jcal^\ell(M,\omega)$ by 
$\Pi_\pi$, which is again a Fredholm map of index $l$.
Consider the subset $\Mod^\circ_\pi(\beta)$
of $\Mod_\pi(\beta)$ consisting of the points of the 
form $\phi=(\pi:\Sig\ra C,\imath_C:C\ra M,J)$ 
such that $\imath_C$ is an embedding. This open subset of $\Mod_\pi(\beta)$
fibers over the subset $\Mod_h^\circ(\alpha)$. 
Abusing the notation, we re-define the bundles $\Bundle$ and
$\Fbundle$ over $\Mod_\pi^\circ(\beta)$ by defining the fibers at a point  
$\phi=(\pi:\Sig\ra C,\imath_C:C\ra M,J)\in \Mod_\pi^\circ(\beta)$ by
\begin{equation}
 \begin{split}
 &\Fbundle_{\phi}=\Gamma^{k,p}(\Sig,\pi^*N_C),\ \ \&\ \ \Bundle_{\phi}=
 \Gamma^{k-1,p}(\Sig,\Omega^{0,1}_\Sig\otimes_J N_C).
 \end{split}
\end{equation}
The bundle $\Bundle$ may be pulled back over $\Fbundle$ using the
projection map  $\Fbundle\ra \Mod_\pi^\circ(\beta)$.
We abuse the notation and will denote this pull back by $\Bundle$ as well.
Then the operator $\Kcal$ defines  a section of $\Bundle\ra \Fbundle$.
\\

For a fixed maximal left ideal 
$\m\in I(\pi)$ of $\R_\Group$ let 
 $\Fbundlebar^\m$  denote the 
\emph{sub-bundle} of $\Fbundle$  consisting of
the points $(\phi;X)$ with $\phi\in\Mod_\pi^\circ(\beta)$ and 
$X\in \Fbundle_\phi$ such that $\m\subset\m_X$.
Let $\Fbundle^\m\subset \Fbundlebar^\m$ denote the 
\emph{subspace} which consists of those   
$(\phi;X)$ with  $\m=\m_X$.
Since $\m$ is maximal, $\Fbundle^\m_\phi$ is the complement of $(\phi;0)$
in $\Fbundlebar^\m_\phi$. Note that the later vector space 
may be identified with $\Gamma^{k,p}_\m(\Sig,\pi^*N_C)$.
Let $\Bundle^\m$ be the sub-bundle of 
$\Bundle|_{\Fbundle^\m}\ra \Fbundle^\m$ consisting of the 
tuples $(\phi;X;\delta)$ with $(\phi;X)\in \Fbundle^\m$ and 
$\delta\in\Bundle_\phi$ so that 
$\m\subset \m_\delta$. This gives a vector bundle 
$\Bundle^\m\ra \Fbundle^\m$ and a section
$\Kcal^\m:\Fbundle^\m\ra \Bundle^m$.
For a point $(\phi;X)\in\Fbundle_\phi^\m$ with
$$\phi=(\pi:\Sig\ra C, \imath_C:C\hookrightarrow M,J)\in
\Mod_\pi^\circ(\beta),\ \ \&
\ \ X\in \Gamma^{k,p}_\m(\Sig,\pi^*N_C)$$
the section $\Kcal^\m$ is defined by
$$\Kcal^\m(\phi;X)=\Kcal^\m_\pi(X)\in\Gamma^{k-1,p}_\m\left(\Sig,
\Omega^{0,1}_\Sig\otimes_J\pi^*N_C\right)=\Bundle_{(\phi;X)}^\m.$$

\begin{prop}\label{prop:intersection-is-transverse}
With the above notation fixed, for any maximal left ideal $\m$ of $\R_\Group$ 
such that the associated 
irreducible representation is faithful and the corresponding matrix algebra is 
either $M_i(\R)$ for $i=1,2,3$ or $M_1(\C)=\C$, the intersection of 
$\Kcal^\m:\Fbundle^\m\ra \Bundle^\m$ with the zero section of the 
vector bundle $\Bundle^\m\ra \Fbundle^\m$ is transverse.
\end{prop}
\begin{proof}
Suppose that $\phi=(\pi,\imath_C,J)$ is a point of 
$\Mod_\pi^\circ(\beta)$ and that $(\phi;X)$
is a point of $\Fbundle^\m$ such that $\Kcal_{\pi}(X)=0$.
We should show that the differential of $\Kcal^\m$, 
projected over the fiber, defines a surjective
 operator
$$d\Kcal^\m:T_{(\phi;X)}\Fbundle^\m\simeq 
T_{\phi}\Mod_\pi(\beta)\oplus \Fbundlebar^\m_{\phi}\ra \Bundle^\m_{\phi}.$$
The tangent space $ T_{\phi}\Mod_\pi(\beta)$ 
has a subspace which consists of the elements $u$ of 
$\Tcal_J=T_J\Jcal^\ell(M,\omega)$
such that $u|_{\Image(\imath_C)}=0$ (and the tangent vector 
is trivial in the direction of $\MapMod$).
Denote this subspace by $\Hcal_J(\imath_C)$. The restriction of $d\Kcal^\m$ 
to $\Hcal_J(\imath_C)\oplus \Fbundlebar^\m_{\phi}$ may be easily computed:
\begin{equation}
 \begin{split}
  &d\Kcal^\m:\Hcal_J(\imath_C)\oplus \Gamma_\m^{k,p}(\Sig,\pi^*N_C)\ra
   \Gamma_\m^{k-1,p}(\Sig,\Omega^{0,1}_\Sig\otimes_J \pi^*N_C)\\
  &d\Kcal^\m(u,Y)=\Kcal^\m(Y)+(\nabla_Xu).d\pi.j.
\end{split}
\end{equation}
The following lemma then completes the proof of this proposition.
\end{proof}

\begin{lem}{\label{lem:nabla-is-surjective}}
With the above notation fixed, 
suppose that $0\neq X\in \mathrm{Ker}(\Kcal_\pi)$ and that  
$\m=\m_X$ is a maximal left ideal of $\R_\Group$ such that the associated 
irreducible representation is faithful and the corresponding 
matrix algebra is either $M_i(\R)$ for $i=1,2,3$
or $M_1(\C)=\C$. Then operator $\nabla_X:\Hcal_J(\imath_C)\ra
 \Coker(\Kcal^{\m}_\pi)$ defined by $\nabla_X(u)=(\nabla_Xu)d\pi.j$
is surjective.
\end{lem}
\begin{proof}
We will present the proof in the cases where the associated matrix algebra of $\m$ 
is $M_3(\R)$ or $M_1(\C)$. The other two cases are in fact easier (and similar).
Let us denote the linear operator $\Kcal^\m_\pi$ by $L$, 
and the normal bundle $N_C$ by $N$
for simplicity. 
Suppose that $\nabla_X$ is not surjective. Identify the
cokernel of $L$ with the kernel of its adjoint $L^*$. 
Choose the non-zero section $\delta$ in 
$\Coker(L)$ so that it is orthogonal to the image of $\nabla_X$.
Since $L^*(\delta)=0$, $\delta$ is not identically zero on any open subset of $\Sig$.
Choose an open ball $U\subset C^\circ=C-\pi(B)$.
Let $p\in U$ be a point in $U$ and let $\{p_\sig\}_{\sig\in\Group}$ be the 
set of pre-images, which are indexed so that $\tau(p_\sig)=p_{\tau\sig}$
for all $\sig,\tau\in\Group$.
We may assume that $\pi^{-1}(U)=\cup_{\sig\in\Group}U_\sig$, with $U_\sig$ a 
ball around $p_\sig$ such that the balls $U_\sig$ are disjoint and 
$\pi:U_\sig\ra U$ is a diffeomorphism. Let $z$ 
denote the local coordinate over $U$ and 
$z_\sig$ denote the corresponding local coordinate over $U_\sig$.
Assume that $\delta$ is non-zero in $U_e$ (where $e$ 
denotes the identity element of $\Group$).
We will denote $U_e$ by $V$ and $z_e$ by $w$.
Fix a bump function $\lambda$ with support inside $U$
 and lift it to $\Sig$, keeping the same
name $\lambda$ for it. By choosing $U$ small enough, we may assume that there is a 
section $u\in \Hcal_J(\imath_C)$, such that $\nabla_X(u)$ is any given section over 
$V$. The section $\lambda.u$ is then supported on $\pi^{-1}(U)$, and  we have
\begin{displaymath}
\begin{split}
0=\langle \nabla_X (\lambda u),\delta \rangle &=\frac{1}{2\sqrt{-1}}
\sum_{\sig\in\Group}\int_{U_\sig} \lambda(z_\sig)\langle \nabla_X (u),
\delta \rangle_{z_\sig}  dz_\sig\wedge d\overline{z}_\sig\\
&=\frac{1}{2\sqrt{-1}}
\int_{V} \lambda(w)\left(\sum_{\sig\in\Group}\langle \nabla_{\sig^*X} (u),
\sig^*\delta \rangle_{w}\right) dw\wedge d\overline{w}.
\end{split}
\end{displaymath}
For the second equality in the above equation we use the fact that 
the restriction of $u$ to the image of $C$ is zero, and the 
function $\lambda$ may thus be brought out of the differentiation.\\
The complex anti-linear $1$-form $\delta$ can locally be written as
$\delta=Zds-JZdt$
where $Z\in \Gamma^{k,q}(\pi^{-1}(U),\pi^* N)$ is a given vector field,
and $z=s+t\sqrt{-1}$ is the local coordinate over $U$.
Thus, for any $p\in C^\circ$ and any linear transformation
$B:N_p\ra N_p$, by letting $\lambda$ converge to 
the delta function supported above $p$
and choosing $u$ so that $\nabla_X(u)(p)$ corresponds 
to $B(X(p))ds-JB(X(p))dt$ we will have
\begin{equation}\label{eq:inner-product}
F_B(p):=\sum_{\pi(q)=p}  \langle BX,Z \rangle_{q}=0.
\end{equation}
If $U$ is  small, $N$ may be trivialized over $U$ as $\C\oplus \C$.
The vector fields $X$ and $Z$ will then have a corresponding decomposition
$X=(X_1,X_2)$ and $Z=(Z_1,Z_2)$. Varying the matrix $B$ in the equation 
$F_B(p)=0$ we conclude that for all $i,j\in\{1,2\}$ and all $p\in C^\circ$
\begin{equation}{\label{eq:local-inner-product}}
F_{ij}(p):=\sum_{\pi(q)=p} (X_i.Z_j)(q)=0\ \&\ \Fbar_{ij}(p):=
\sum_{\pi(q)=p} (X_i.\overline{Z}_j)(q)=0.
\end{equation}
Let us first assume that the matrix algebra associated with $\m$ is $M_1(\C)=\C$.
Since the corresponding representation is faithful, we have 
$$\Group=\Z/n\Z=\langle \sig\rangle,\ \ 
\mathrm{with }\ \rho(\sig)=\zeta\in S^1-\{\pm 1\}\subset \C^*.$$
Here $\zeta$ denotes a primitive $n$-th root of unity. We may thus find a second 
section $Y\in \Gamma^{k,p}(\Sig,\pi^*N)$ such that
\begin{equation}
 \label{eq:action-on-X}
\left(\begin{array}
{c}(\sig^m)^*X\\ (\sig^m)^*Y 
\end{array}\right)
=\left(\begin{array}
        {cc} \Real(\zeta^m)&\Image(\zeta^m)\\ -\Image(\zeta^m)&\Real(\zeta^m)
       \end{array}\right)
\left(\begin{array}
       {c} X\\Y
      \end{array}\right).
\end{equation}
Similarly, there is a section $\epsilon\in\Gamma^{k,q}(\pi^{-1}(U),
\Omega^{0,1}_\Sig\otimes_J\pi^*N)$
such that 
\begin{equation}
 \label{eq:action-on-d}
\begin{split}
&\epsilon=Wds-JWdt, \ \ W\in\Gamma^{k,q}(\pi^{-1}(U),\pi^*N)\\
&\left(\begin{array}
{c}(\sig^m)^*Z\\ (\sig^m)^*W 
\end{array}\right)
=\left(\begin{array}
        {cc} \Real(\zeta^m)&\Image(\zeta^m)\\ -\Image(\zeta^m)&\Real(\zeta^m)
       \end{array}\right)
\left(\begin{array}
       {c} Z\\W
      \end{array}\right).
\end{split}
\end{equation}
Rewriting equation~\ref{eq:local-inner-product} 
using equations~\ref{eq:action-on-X} and~\ref{eq:action-on-d}
we obtain
\begin{equation}
\label{eq:local-inner-product-2} 
 \begin{split}
 0&=\sum_{m=0}^{n-1}\Big(\Real(\zeta^m)^2X_iZ_j+\Image(\zeta^m)^2Y_iW_j\\
&\ \ \ \ \ \ +\Image(\zeta^m)\Real(\zeta^m)(X_iW_j+Y_iZ_j)\Big)_q\\
&= \left(\sum_{m=0}^{n-1}\Real(\zeta^m)^2\right)
(X_iZ_j)_q+\left(\sum_{m=0}^{n-1}\Image(\zeta^m)^2\right)(Y_iW_j)_q\\
&\ \ \ \ \ \ +\left(\sum_{m=0}^{n-1}\Image(\zeta^m)\Real(\zeta^m)\right)
(X_iW_j+Y_iZ_j)_q\\
&=\left(\sum_{m=0}^{n-1}\Real(\zeta^m)^2\right)(X_iZ_j+Y_iW_j)_q,
\ \ \ \ \ \ \ \ \ \ \ 
\ \ \ \ \ \ \ \ \ \ \ \ \ \forall\ \ q\in V. 
 \end{split}
\end{equation}
Here the last equality follows since from $n\neq 2$ we have 
$\sum_{m=0}^{n-1} \zeta^{2m}=0$ and
$$\zeta^{2m}=\left(\Real(\zeta^m)^2-\Image(\zeta^m)^2\right)
-2\sqrt{-1}\left(\Real(\zeta^m)\Image(\zeta^m)\right).$$
Since $\sum_{m=0}^{n-1}\Real(\zeta^m)^2$ is a positive real number, 
we conclude that 
$X_iZ_j+Y_iW_j$ is identically zero on $V$.
Similarly, from the second equality in equation~\ref{eq:local-inner-product} 
we may conclude that
$$\left(X_i\ov{Z_j}+Y_i\ov{W_j}\right)\big|_{V}\cong 0.$$
These two relations imply that 
\begin{equation}\label{eq:local-inner-product-3}
 \left(\Real(Z_j).X+\Real(W_j)Y\right)(w)=0,\ \ \forall\ w\in V.
\end{equation}
Since $X$ and $Y$ are both in the kernel of $L$, and the zeros of both of them are 
isolated, equation~\ref{eq:local-inner-product-3}
implies that $\dbar_w\left(\frac{\Real(W_j)}{\Real(Z_j)}\right)=0$ on $V$, and hence 
$\Real(W_j)/\Real(Z_j)$ is constant on $V$ 
(as a real valued holomorphic function). This means that $Y$ is a constant multiple of 
$X$ over $V$, and hence on all of $\Sig$,
unless $\Real(Z_j)=\Real(W_j)=0$. Since the former can not happen, we should have 
$\delta=0$ on $V$, which is also a contradiction.
This completes the proof when the matrix algebra 
is $\C$.\\
\\
Let us now assume that the matrix algebra 
associated with $\m$ is $M_3(\R)$. Then associated with any 
$\sig\in\Group$ is an orthogonal matrix $A_\sig\in O_3(\R)\subset M_3(\R)$.
Furthermore, assuming $\m$ corresponds to the matrices with 
zero in the first column, we may write 
\begin{equation}
 \label{eq:3-term-decomposition}
 \begin{split}
  &\alpha(X)=X^1\alpha_1+X^2\alpha_2+X^3\alpha_3,\ \ \&\\ 
&\alpha(Z)=Z^1\alpha_1+Z^2\alpha_2+Z^3\alpha_3,\\ 
\end{split}
\end{equation}
where we have 
\begin{displaymath}
 \begin{split}
  &\alpha_1=\left(\begin{array}
                                   {ccc}1&0&0\\ 0&0&0\\ 0&0&0
                                  \end{array}\right),
\                 \alpha_2=\left(\begin{array}
                                   {ccc}0&1&0\\ 0&0&0\\ 0&0&0
                                  \end{array}\right),
\                 \alpha_3=\left(\begin{array}
                                   {ccc}0&0&1\\ 0&0&0\\ 0&0&0
                                  \end{array}\right),\\
&X^i\in\Gamma^{k,p}(\Sig,\pi^*N)\ \ \ \&\ \ \ 
Z^j\in\Gamma^{k,q}(\pi^{-1}(U),\pi^*N),\ \ \ \forall\ i,j\in\{1,2,3\}.
 \end{split}
\end{displaymath}

From the presentations of equation~\ref{eq:3-term-decomposition}, setting
$\beta(X)=(X^1, X^2, X^3)^t$ and $\beta(Z)=(Z^1,Z^2,Z^3)^t$ we  have 
$$\sig^*(\beta(X))=A_\sig.\beta(X),\ \ \ \&\ \ \ 
\sig^*(\beta(Z))=A_\sig.\beta(Z).$$
There is a fixed vector $v=(v_1,v_2,v_3)\in\R^3$ such that $X=v.\beta(X)$ 
and $Z=v.\beta(Z)$.
If $X_i$ and $Z_i$ denote the $i$-th components of $X$ and $Z$ in a trivialization of 
$N$ over $U$ as $\R^4\times U$, for a fixed point $q\in V$ which satisfies 
$\pi(q)=p$, equation~\ref{eq:local-inner-product} reads as
\begin{equation}
 \label{eq:local-inner-product-4}
 \begin{split}
  0&=\sum_{\pi(z)=p}Z_i(z)X(z)=\sum_{\sig\in\Group}Z_i(\sig(q))X(\sig(q))\\
&=\beta(X)^{t}(q).\left(\sum_{\sig\in\Group}A_\sig^t.v^t.v.
A_\sig\right).\beta(Z_i)(q)=:
\beta(X)^{t}(q).B.\beta(Z_i)(q).
 \end{split}
\end{equation}
For the $3\times 3$ matrix $B$ and any choice of $\tau\in\Group$ we have
\begin{displaymath}
\begin{split}
 &A_{\tau}^{-1}BA_\tau=A_\tau^tBA_\tau=\sum_{\sig\in\Group}
 A_{\sig\tau}^t.v^t.v.A_{\sig\tau}=B\ 
 \Rightarrow\ \  BA_\tau=A_\tau B,\ \ \forall\ \tau\in\Group.
\end{split}
\end{displaymath}
Since the matrices $\{A_\tau\}_{\tau\in\Group}$ generate the algebra $M_3(\R)$,
$B$ should be a multiple of identity. On the other hand, 
the trace of $B$ may be computed via
\begin{equation}
 \label{eq:trace-of-B}
\begin{split}
&\trace(B)=\sum_{\sig\in\Group}\trace(A_\sig^t.v^t.v.A_\sig)
=\sum_{\sig\in\Group}\trace(v.A_\sig.A_\sig^t.v^t)=3|\Group|\|v\|^2\neq 0\\
&\Rightarrow\ \ \ \ B=\|v\|^2|\Group| .I_{3\times 3}\neq 0.
\end{split}
\end{equation}
Combining equations~\ref{eq:local-inner-product-4} and ~\ref{eq:trace-of-B}
it follows that $\beta(X)\beta(Z_i)^t$ is identically zero 
over $\pi^{-1}(U)$ for $i=1,2,3,4$.\\
The three vectors $X^1(q), X^2(q)$ and $X^3(q)$ will thus 
linearly depend on each other (over the 
point $q$). If $r(q)$ denotes the rank of the vector 
space spanned by these three vectors we will have $r(q)\in\{1,2\}$
for a generic choice of $q$. If $r(q)=2$, then for points in an 
open neighbourhood of $q$ the same will be true. 
For any point $z$ in this open neighbourhood, the four vectors 
$$\beta(Z_i)(z)=(Z_i^1(z),Z_i^2(z),Z_i^3(z))\in\R^3, \ \ \ i=1,2,3,4$$  
are thus multiples of one-another, and $Z^1$ is thus a real multiple of $Z^2$
over this open neighbourhood, i.e. $Z^1(z)=\lambda(z).Z^2(z)$ for 
some real valued function $\lambda$.
Since $Z^1$ and $Z^2$ satisfy perturbed Cauchy-Riemann 
equations, this means that $\partial_z(\lambda)=0$,
implying that $\lambda$ is constant (since it is real valued). 
Repeating this argument for the other pairs,
it is implied that there is a constant vector 
$0\neq w=(w_1,w_2,w_3)\in\R^3$ such that $w.\beta(X)=0$ over 
a small open neighbourhood on $\Sig$, and hence on all of $\Sig$. We thus have
\begin{displaymath}
 \begin{split}
  &0=\sig^*\left(w.\beta(X)\right)=(w.A_\sig).\beta(X)\ \ \ \forall\ \sig\in\Group\\
  &\Rightarrow\ \ \ \beta(X)=0, \ \ \text{which implies }\ X=0.
 \end{split}
\end{displaymath}
From this contradiction, we should have $r(q)=1$ for a generic point $q$ on $\Sig$.
If $r\leq 1$ in a neighbourhood of $q$, it is implied that  
$X^2(z)=\lambda(z) X^1(z)$ for $z$ near $q$ and for a real valued function 
$\lambda$. Again, since $X^1$ and $X^2$ satisfy 
perturbed Cauchy-Riemann equations, this means that $\dbar_z(\lambda)=0$,
implying that $\lambda$ is constant near $q$ (since it is real valued).
The equation $X^2(z)=\lambda X^1(z)$ (with $\lambda$  a real constant) thus 
extends to all of $\Sig$, and the same argument as above may be repeated for 
$w=(-\lambda,1,0)$. The resulting contradiction then completes the proof.
\end{proof}

\begin{thm}{\label{thm:main}}
 For any $J$ in  a subset
$\Jcal^{\ell}_{s.r}(M,\omega)\subset \Jcal^\ell_{reg}(M,\omega)$
which is of second category as a subset of $\Jcal^\ell(M,\omega)$,
any genus $h\geq 0$, and any homology class $\alpha\in\Ht_2(M,\Z)$ which satisfies
$c_1(M).\alpha=0$, all the curves in the moduli space $\Mod_h(\alpha,J)$ 
are smooth  $4$-rigid embeddings.
\end{thm}
\begin{proof}
Let $\alpha\in\Ht_2(M,\Z)$ be a class which satisfies 
$c_1(M).\alpha=0$, and $g\geq h\geq 0$ be fixed integers.
Let $\pi:\Sig\ra \Sig/\Group=C$ denote the topological type
of a branched covering map coming from a group action on a surface $\Sig$ of genus
$g$ and with quotient a surface $C$ of genus $h$. Furthermore, assume that 
$\Group$ is a subgroup of $S_4$. Let $\beta=|\Group|.\alpha\in\Ht_2(M,\Z)$.
Fix $\m\in I(\pi)$ and assume that the corresponding representation is faithful. Let
 $\Nod_\pi^\m(\beta)$ be the zero locus of the section $\Kcal^\m_\pi$, 
 which is fibered
over $\Jcal^\ell(M,\omega)$.
By proposition~\ref{prop:intersection-is-transverse} each one of the moduli
spaces $\Nod_\pi^\m(\beta)$ is a smooth infinite dimensional manifold.
The set of regular values of the Fredholm projection map
$$\Pi_\pi^\m(\beta):\Nod_\pi^\m(\beta)\ra \Jcal^\ell(M,\omega)$$ 
is of second category. Denote the set of
regular values for the projection map 
$\Pi_\pi^\m(\beta)$ by $\Jcal_{reg}^\ell(\pi,\m,\alpha)$ and set
\begin{equation}\label{eq:regular-ACSs}
\Jcal^\ell_{n}(M,\omega)
=\Jcal^{\ell}_{reg}(M,\omega)\bigcap 
\left(\bigcap_{\pi,\m}\bigcap_{\substack{\alpha\in\Ht_2(M,\Z)\\ c_1(M).
\alpha=0}}\Jcal_{reg}^\ell(\pi,\m,\alpha)\right).
\end{equation}
Here, the first intersection inside the parentheses is over all 
the topological types of the 
branched covering maps 
$\pi$ coming from the action of a subgroup of $S_n$ and all $\m\in I(\pi)$ which 
correspond to faithful representations. When $n=4$, the corresponding 
matrix algebra of $\m$ is 
thus one of $M_i(\R),\ i=1,2,3$ or $M_1(\C)=\C$.\\

For an almost complex structure $J\in \Jcal^{\ell}_4(M,\omega)$ 
and any homology class
$\alpha\in\Ht_2(M,\Z)$ satisfying $c_1(M).\alpha=0$, $\Mod_h(\alpha,J)$
 is a zero dimensional manifold
since $\Jcal_4(M,\omega)\subset \Jcal^\ell_{reg}(M,\omega)$. If an embedded curve 
$C\in\Mod_h(\alpha,J)$ is not $4$-rigid, we may construct a branched covering map 
$\pi:\Sig\ra \Sig/\Group=C$ for a 
subgroup $\Group<S_4$, so that the kernel of $\Kcal_\pi^\m$ is 
non-trivial for some $\m\in I(\pi)$.
Furthermore, we may assume that the representation associated with $\m$ is faithful
(all these reductions were discussed in section~\ref{sec:reductions}).
This non-trivial kernel gives a non-empty subset of $\Nod_\pi^\m(\beta,J)$.
The moduli space $\Nod_\pi^\m(\beta,J)$ is a
smooth manifold of dimension $\Index\left(\Pi_\pi^\m(\beta)\right)$, by regularity
of the almost complex structure $J$. This later index 
 is equal to the sum of the indices of the projection
map from $\Mod_\pi(\beta)$ to $\Jcal^\ell(M,\omega)$ 
and the projection map from $\Nod_\pi^\m(\beta)$ to 
$\Mod_\pi(\beta)$. This later projection map has 
the same index as the operator $\Kcal^\m$.
In the following section, we will estimate $\Ind(\Kcal^\m)$ and will prove
proposition~\ref{prop:index-of-Ktild}.
\begin{prop}\label{prop:index-of-Ktild}
If $\m\in I(\pi)$ corresponds to a faithful representation,  
the index $\Ind(\Kcal^\m)$ is at most $-l$, where
$l$ is the number of points in the image $\pi(B)$ of the branched locus of $\pi$.
\end{prop}
Assuming proposition~\ref{prop:index-of-Ktild},
the index of $\Pi_\pi^\m(\beta)$ is equal to
$l+\Ind(\Kcal^\m)$, which is  at most $0$.
If this manifold is non-empty, and $(\phi;Y)\in \Nod_\pi^\m(\beta,J)$, 
then for any $0\neq\lambda\in \Rbb$, $(\phi;\lambda Y)$ is also in
 $\Nod_\pi^\m(\beta,J)$. This contradicts the fact  that this moduli
space is at most zero dimensional. The contradiction proves that every curve 
$C\in\Mod_h(\alpha,J)$ is $4$-rigid.
\end{proof}

\section{Index computation}\label{sec:index}
Fix a holomorphic map $\pi:\Sig\ra C$ of degree $n>1$ which is 
determined by the action of
a group $\Group$ on the Riemann surface $\Sig$ (so that $C=\Sig/\Group$).
Let $L\ra C$ be a holomorphic vector bundle over $C$ and 
$\pi^*L\ra \Sig$ be its pull-back 
over $\Sig$. The group ring $\R_\Group$ acts on the 
cohomology groups $\Ht^i(\Sig,\pi^*L)$. 
For $\eta\in\Ht^i(\Sig,\pi^*L)$ we may set $\m_\eta$ to be the left ideal of 
$\R_\Group$ 
consisting of the elements $\afrak\in\R_\Group$ such that $\afrak.\eta=0$.
 For any  $\m\in I(\pi)$ and any vector bundle of the form $\pi^*L$ let
\begin{displaymath}
\Ht^i_{\mathfrak{m}}(\pi^*L):=\big\{\eta \in \Ht^i(\Sigma,\pi^*L) \ \big | \
\mathfrak{m}\subset \m_\eta \big \},\hspace{1cm} i=0,1.
\end{displaymath}
Define
$h^i_{\mathfrak{m}}(\pi^*L)=\text{dim}_{\mathbb{C}}
(\Ht^i_{\mathfrak{m}}(\pi^*L)),\ i=0,1$, and let
$\chi_{\mathfrak{m}}(\pi^*L)=
h^0_{\mathfrak{m}}(\pi^*L)-h^1_{\mathfrak{m}}(\pi^*L)$.
For any $\m\in I(\pi)$ the Fredholm map $\Kcal^\m=\Kcal_\pi^\m$ 
is a zero order perturbation of the operator  
$$\dbar:\Gamma^{k,p}_\m(\Sig,\pi^*N_C)\ra 
\Gamma^{k-1,p}_{\m}(\Sig,\Omega^{0,1}_\Sig\otimes_{J}\pi^*N_C),$$
i.e. the operator obtained when $J$ is the integrable complex structure on a 
neighbourhood of the zero section in $N_C$.
The index of $\Kcal^\m$ may thus be computed 
by a computation for $\dbar$.
In this case (i.e. for $\Kcal^\m=\dbar$) we have
$$\Ker(\Kcal^\m)=\Ht^0_\m(\Sig,\pi^*N_C),\ \ \&\ \ 
\Coker(\Kcal_\m)=\Ht^1_\m(\Sig,\pi^*N_C).$$
Since $N_C$ may be deformed to $\Ocal_C\oplus K_C$, 
where $K_C$ denotes the canonical bundle of $C$ 
(see \cite{R-Jim} and \cite{P}), we have 
$$\Index(\Kcal^\m)=\chi_\m(\pi^*K_C)+\chi_\m(\pi^*\Ocal_C).$$
Let $B$ denote the branching divisor of the map 
$\pi$ and $\Ocal_B$ denote the corresponding sheaf 
over $\Sig$ with support in $B$. The group ring $\R_\Group$ 
will also act on $\Ht^0(\Sig,\Ocal_B)$ and 
it thus makes sense to talk about $\Ht^0_\m(\Ocal_B)$ and $h^0_\m(\Ocal_B)$.
\begin{lem}\label{lem:index}
For  any $\mathfrak{m}\in I(\pi)$ we have  $\Ind(\Kcal^\m)=-2h^0_\m(\Ocal_B)$.
\end{lem}
\begin{proof} 
Since $K_\Sig$ is invariant under the action of $\Group$, we may apply Serre duality
to obtain $\chi_\m(\pi^*K_C)=-\chi_\m(K_\Sig-\pi^*K_C)=-\chi_\m(B)$. The
line bundle associated with $B$  sits in a short exact sequence
\begin{displaymath}
 0\lra \Ocal_\Sig\lra B\lra \Ocal_B\lra 0.
\end{displaymath}
Considering the cohomology long exact sequence
associated with this short sequence, and restricting attention to the 
sections trivialized by the maximal left ideal $\m$, 
we obtain $\chi_\m(B)=\chi_\m(\Ocal_\Sig)+h^0_\m(\Ocal_B)$.
Since $\pi^*\Ocal_C=\Ocal_\Sig$, the proof is complete.
\end{proof}
\begin{proof}(of proposition~\ref{prop:index-of-Ktild})
By the above lemma, the computation of the index 
of the operator $\Kcal^\m$ is reduced to 
the computation of $h^0_\m(\Ocal_B)$. 
Let $x\in  C$ be a point in the image of the branched locus of $\pi$.
The monodromy around $x$ gives a partition of $n$ as 
$n=d_1+d_2+...+d_\ell$  with
$d_1\geq d_2\geq ...\geq d_\ell$. Associated with each $d_i$ is a point 
$y_i$ in $\pi^{-1}(x)$ which has
multiplicity $d_i$ and is a branched point if $d_i>1$. 
Let $B_x$ be the divisor $(d_1-1)y_1+...+(d_\ell-1)y_\ell$.
Since $\pi$ comes from the action of a group, we should 
have $d_1=...=d_m=d>1$, and that the degree of $B_x$ is 
$m(d-1)$. Set $\Bcal_x=\Ocal_{B_x}$ and note that $B=\sum_x B_x$ and 
$\Ocal_B=\sum_x\Bcal_x$.  From $\Ht^0_\m(\Ocal_B)
=\oplus_x\Ht^0_\m(\Bcal_{x})$, it suffices to show $h^0_\m(\Bcal_{x})\geq 1$.\\
The elements of $\Ht^0(\Bcal_x)$ are the germs
$\phi=\sum_{i=1}^{\ell}\sum_{j=1}^{d-1}{A_{i,j}}/{z_i^j}$,
where $z_i$ is the pull-back of the local coordinate 
$z$ around $x$ to a neighborhood of $y_i$.
If $\zeta=\exp(2\pi\sqrt{-1}/d)$ is the primitive $d$-th root of unity, 
there is a unique
$\tau\in\Group$ such that $\tau(z_1)=\zeta z_1$ is satisfied near $y_1$.
We may also assume, after re-parametrization, that there are elements 
$\sig_1,...,\sig_m\in\Group$ such that $\sig_j(z_1)=z_j$ is satisfied for $j=1,...,m$.
Note that $\sig_1$ would be the identity element of $\Group$. The group $\Group$
is then identified as
$$\Group=\Big\{\sig(i,j):=\sig_i\tau^j\ \big|\ i=1,...,m,\ j=1,...,d\Big\}.$$  
For $(i,j)\in\{1,...,m\}\times\{1,...,d\}$ we define
$\xfrak_{i,j}:=\sum_{p=1}^{d-1}\zeta^{pj}z_i^{-p}\in\Ht^0(\Bcal_x)$.
For $\sig\in\Group$, $\sig(z_i)=\zeta^az_{\sig(i)}$, where $a\in\{1,...,d\}$ and 
$\sig(i)\in\{1,...,m\}$ are integers which depend on $\sig$ and $i$. Thus
\begin{displaymath}
\sig(\xfrak_{i,j})=\sum_{p=1}^{d-1}\zeta^{pj}\left(\zeta^a.\zeta_{\sig(i)}\right)^{-p}
=\sum_{p=1}^{d-1}\zeta^{p(j-a)}z_{\sig(i)}^{-p}=\xfrak_{\sig(i),j-a}.
\end{displaymath}
Thus $\Group$ permutes the elements $\xfrak_{i,j}$ among themselves. If we denote 
$\sig(\xfrak_{1,d})$ by $\xfrak_\sig$ we will have 
$$\sig_1(\xfrak_{\sig_2})=\xfrak_{\sig_1\sig_2},\ \ \ \forall\ 
\sig_1,\sig_2\in\Group.$$
There are precisely $m$ relations among the sections 
$\xfrak_\sig$, indexed by $i=1,...,m$:
\begin{displaymath}
\begin{split}
 \sum_{j=1}^d\xfrak_{\sig(i,j)}&=
\sig_i\left(\sum_{j=1}^{d}\sum_{p=1}^{d-1}\tau^j(z_1)^{-p}\right)=
\sig_i\left(\sum_{p=1}^{d-1}z_1^{-p}\left(\sum_{j=1}^{d}\zeta^{-jp}\right)\right)=0.
\end{split}
\end{displaymath}
The sections $\{\xfrak_\sig\}_{\sig\in\Group}$ generate 
$\Ht^0(\Bcal_x)$. An element $\phi\in\Ht^0(\Bcal_x)$ of the form 
$\phi=\sum_{\sig\in\Group}a_\sig\xfrak_\sig$ is in $\Ht^0_\m(\Bcal_x)$ 
if and only if
$$\sum_{\sig\in\Group}b_\sig a_\sig\xfrak_\sig=0,\ \ \ \ 
\forall \ \ \bfrak=\sum_{\sig\in\Group}b_\sig\sig^{-1}\in\m.$$
This means that $\Ht^0_\m(\Bcal_x)$ is non-trivial if and 
only if the sum of the subspace
$\m$ of the vector space $\R_\Group$ with the subspace 
$$\left\langle\alpha(i):=\sum_{j=1}^d\sig(i,j)^{-1}\ \Big|\
 i=1,...,m\right\rangle_\R$$
is not all of $\R_\Group$. Note that this sum should be considered 
as the sum of sub-vector-spaces, and not ideals.\\ 
Let us assume that the vector space $\m\in\R_\Group$ is of 
codimension $\ell$. Thus the subspace $\R_\Group/\m$ 
generated by $\sig_1^{-1},...,\sig_m^{-1}$ is at most of dimension 
$\ell$, and we may assume that it 
is generated by $\sig_1^{-1},...,\sig_\ell^{-1}$. Thus there are real numbers 
$c_{i,k}\in\R$, $i=1,...,m$ and $k=1,...,\ell$ such that 
$\sig_i^{-1}=\sum_{k=1}^\ell c_{i,k}\sig_k^{-1}$ for $i=1,...,m$. We thus have
\begin{displaymath}
\alpha(i)=\sum_{k=1}^\ell c_{i,k}\alpha(k),\ \ (\mathrm{modulo}\ \m),\ \ \ 
\forall\ \ i=1,...,m,
\end{displaymath}
which implies that $\Ht_\m^0(\Bcal_x)\neq 0$ if and only if (again as vector spaces)
$$\m+\langle \alpha(i)\ |\ i=1,...,\ell\rangle_\R\neq \R_\Group.$$
It thus suffices to show that there is some non-zero 
$a=(a_1,...,a_\ell)\in\R^\ell$ such that
$$\sum_{i=1}^\ell a_i\alpha(i)=\left(\sum_{j=0}^{d-1}\tau^j\right)
\left(\sum_{i=1}^{\ell}a_i\sig_{i}^{-1}\right)\in\m.$$
This criteria should be investigated over the distinguished component
of the maximal left ideal $\m$ in the decomposition of
 equation~\ref{eq:decomposition-R_G}.
Over this component, we have an isomorphism with a matrix algebra of the 
form $M_k(R)$ for an integer $k=\ell,\ell/2$ or $\ell/4$ 
 depending on whether $R=\R,\C$ or $\Hbb$, respectively.
The ideal $\m$ determines a left ideal of $M_k(R)$, which may be 
specified by a vector 
$v_\m\in R^k$ in the sense that
$$\m\cap M_k(R)=\{A\in M_k(R)\ |\ A.v_\m=0\}.$$
Let $S_1,...,S_\ell$ denote the matrices in $M_k(R)$ 
which correspond to $\sig_1,...,\sig_\ell$. Suppose that the 
criteria is not satisfied. Thus $e_i=S_iv_\m$, for $i=1,...,\ell$ generate
$R^k$ over $\R$. Let $A$ be the matrix in $M_k(R)$ which corresponds to 
$(1+\tau+...+\tau^{d-1})/d$. The above criteria would be satisfied 
if $A(a_1e_1+...+a_\ell e_\ell)=0$. In Other words, we only need to show that 
the kernel of $A$ is non-trivial. Suppose otherwise that $A$ is non-singular.
Since $A^2=A$ we should have $A=I_{M_k(R)}$ is the identity matrix in $M_k(R)$.
Let $T$ denote the matrix which corresponds to $\tau$. Then we have
$$\tau.\frac{1+\tau+...+\tau^{d-1}}{d}=\frac{1+\tau+...+\tau^{d-1}}{d}\ 
\Rightarrow\ TA=A\ \Rightarrow \ T=I_{M_k(R)}.$$
In other words, $\tau$ is in the kernel of the 
representation corresponding to $\m$, which is a 
contradiction (with the assumption that this representation is faithful).
The contradiction completes the proof.
\end{proof}
\begin{remark}
Note that the index computation does not make any use of the 
assumptions on the size of the irreducible real representation, 
and is thus completely general.
\end{remark}
\section{Elliptic regularity and smooth complex structures}\label{sec:regularity}
So far, we have only considered $C^\ell$ almost complex 
structures and $n$-rigidity for embedded curves, and 
sections over them which are of class $W^{k,p}$, where $k\leq \ell$.
Note that if $J$ is an almost complex structure in 
$\Jcal^\ell(M,\omega)$ and $f:\Sig\ra M$ is a $J$-holomorphic curve
of class $W^{1,p}$ for some $p>2$, then it will be of class $W^{\ell+1,p}$ 
and there is a constant 
$c(J,\ell)$ such that $\|f\|_{W^{\ell+1,p}}\leq c(J,\ell)\|f\|_{W^{1,p}}$.\\

 Fix $\ell\in\Z^{>0}\cup\{\infty\}$,  the integers $n,K>0$, 
and an almost complex structure $J\in\Jcal^\ell(M,\omega)$. 
Fix an identification of $\Ht_2(M,\Q)$ with $\Q^{b_2(M)}$ once for all,
and for a homology class $\beta\in\Ht_2(M,\Z)$ let $\|\beta\|$ 
denote its Euclidean norm under this identification. Let $|\beta|$
denote the divisibility of $\beta$.
\begin{defn}
 A $J$-holomorphic map $f:C\ra M$ 
of class $W^{k,p}$ with (possibly disconnected and nodal) domain $C$ 
of total genus $h\leq K$ which is somewhere injective on any 
irreducible component of $C$ is 
called a {\emph{critical}} $(J,n,K)$ {\emph{curve}} if $f(C)$ is 
connected, $\|df\|_{\infty}<K$, $\|f_*[C]\|<K$, $|f_*[C]|\leq n$ 
and either of the following happens:\\
1. The domain $C$ is nodal (i.e. not smooth)\\
2. $C$ is smooth, but there are $x,y\in C$ with $x\neq y$ and $f(x)=f(y)$.\\
3. $C$ is smooth and $f$ is one-to-one, but there is some $x\in C$ with $df_x=0$.\\
4. $C$ is smooth and $f$ is an embedding, but the image of $f$ is not $(n,K)$-rigid.
\\
The $J$-holomorphic map $f:C\ra M$ is called a 
{\emph{critical $(J,n)$ curve}} if it is a critical 
$(J,n,K)$ curve for some $K>0$.
\end{defn}
Thus $\Jcal^\ell_n(M,\omega)$ is included in the subset of 
$\Jcal^\ell(M,\omega)$ consisting of the almost complex structures $J$
so that there are no critical $(J,n)$ curves in $M$. Abusing the notation, we will use
$\Jcal^\ell_n(M,\omega)$ to denote this later (bigger) 
subset of $\Jcal^\ell(M,\omega)$, which is of second category by 
theorem~\ref{thm:main} in $\Jcal^\ell(M,\omega)$, provided 
that $\ell\neq \infty$ and $n\leq 4$. 
Let $\Jcal_n^\ell(K)\subset\Jcal^\ell(M,\omega)$ be the 
subset of almost complex structures $J$ of class $C^\ell$ such that
there are no critical $(J,n,K)$ curves $f:C\ra M$ of class $W^{\ell+1,2}$. 
We will simply denote $\Jcal^\infty_n(K)$ by $\Jcal_n(K)$.
It is then clear that $\Jcal_n^\infty(M,\omega)
=\cap_{K\in \Z^{>0}}\Jcal_n(K)$.\\

\begin{lem}\label{lem:open}
 For any positive real number $K>0$, $\Jcal_n(K)$ is 
 open in $\Jcal^\infty(M,\omega)$ with respect to $C^\infty$ topology.
\end{lem}
\begin{proof}
 Suppose otherwise that there is a sequence $J_i$ of 
 smooth almost complex structures outside $\Jcal_n(K)$ 
which converges to a point $J\in\Jcal_n(K)$. We may 
thus find a sequence $f_i:\Sig_i\ra M$ of critical $(J_i,n,K)$ curves satisfying
$$\|df_i\|_\infty\leq K,\ \ \|(f_i)_*[\Sig]\|\leq K,\ \ \&\ \ |(f_i)_*[\Sig]|\leq n.$$
Furthermore, if $f_i$ is an embedding, it is not $(n,K)$-rigid.
After passing to a subsequence, we may assume that all $f_i$ 
represent the same homology class $\beta$ with $|\beta|\leq n$, and 
that their domains $\Sig_i$ have the same genus $g\leq K$. 
The sequence $\{f_i\}_i$ will thus converge to a $J$-holomorphic curve 
$f:\Sig\ra M$ which satisfies
$$\|df\|_\infty\leq K,\ \ \|f_*[\Sig]\|=\|\beta\|\leq K,\ \ 
\&\ \ |f_*[\Sig]|=|\beta|\leq n.$$
Furthermore, the genus of $\Sig$ is at most $K$. 
The map $f$ decomposes as $f=\imath_C\circ \pi$ where 
$\pi:\Sig\ra C$ is a branched covering map and $\imath_C:C\ra M$
 is $J$-holomorphic and 
somewhere injective on all its components. Since $J\in\Jcal_n(K)$, 
$C$ is smooth, $\imath_C$ is an embedding and the 
image of $C$ under $\imath_C$ is $(n,K)$-rigid. Note also that $\deg(\pi)\leq n$.\\
We may identify a neighbourhood of $C$ with a neighbourhood 
of the zero section in its normal bundle $N_C$ in $M$.
For $i$ sufficiently large, the image of $f_i$ will be in 
this neighbourhood and we may thus regard 
the images of the critical $(J_i,n,K)$ curves $f_i:\Sig_i\ra M$ as 
multi-sections of this normal bundle.
As in section~\ref{sec:moduli-space} we may re-scale these 
images so that their supremum norm is equal to $1$.
Similar to our discussion in section~\ref{sec:moduli-space} this 
convergence gives a section
$$(X,\jmath)\in\Gamma^\infty(\Sig,\pi^*N_C)\oplus \Hom_J^{0,1}(T_C,N_C),$$
which satisfies the following equation in $\Gamma^\infty(\Sig,
\Omega^{0,1}_\Sig\otimes_J\pi^*N_C)$:
\begin{displaymath}
 \begin{split}
 \Kcal_\pi(X)+\jmath\circ d\pi\circ j_\Sig=0.
\end{split}
\end{displaymath}
Here $\Hom_J^{0,1}(T_C,N_C)$ is the space of bundle 
homomorphisms $\jmath$ from the 
tangent space $T_C$ of $C$ to the normal bundle $N_C$
which anti-commute with the complex structure, i.e. 
$J\circ \jmath+\jmath\circ j_C=0$.\\
Since $C$ is $(n,K)$-rigid, for any branched covering map $\varsigma:S\ra C$
such that $\deg(\varsigma)\leq n$ and the genus $g_S$ of $S$ satisfies 
$g_S\leq K+g_C$ the operator $\Kcal_\varsigma$ is injective. In particular, 
for $\varsigma=Id_C:C\ra C$, since the index of $\Kcal_\varsigma$ is zero, 
this operator is an isomorphism. Thus there is some $Y\in\Gamma(C,N_C)$
such that $\Kcal_\varsigma(Y)+\jmath\circ Id_{T_C}\circ j_C=0$.
This means that $X-\pi^*Y\in \Ker(\Kcal_\pi)$. By $(n,K)$-rigidity of 
$C$ we have $X=\pi^*Y$.\\
Since $\imath_C:C\ra M$ is $(n,K)$-rigid 
(and hence $1$-rigid), there are $J_i$-holomorphic
curves $C_i$ close to $C$, which converge to $C$ as $i$ goes to infinity. 
We may then regard both $C_i$ and $\Sig_i$ as sections and 
multi-sections of $N_C$. For large 
values of $i$, the projection map $\imath_i$ from $C_i$ to $C$ is a 
diffeomorphism, and the 
projection map from $\Sig_i$ to $C$ may thus be composed with 
$\imath_i^{-1}$ to give 
the projection maps $\pi_i:\Sig_i\ra C_i$. The map $\pi_i$ may be used to define a
new complex structure $j_i$ on $\Sig_i$ which makes $\pi_i$ 
holomorphic. Next, the bundle 
$N=N_C$ may be pulled back over $C_i$ using $\imath_i$. 
Finally, the map $f_i:\Sig_i\ra M$ may be 
regarded as a section $X_i\in \Gamma^\infty(\Sig_i,\pi_i^*N)$. 
The section $X_i$ describes 
the distance of $\Sig_i$ from $C_i$. Since $X$ is of the form $\pi^*Y$, we have
$$\lim_{i\ra \infty}\frac{\|X_i\|_\infty}{\|J_i-J\|_{C^0}}=0.$$
Since $\Sig_i$ is $J_i$-holomorphic we have
\begin{displaymath}
\begin{split}
 \nabla X_i+J\nabla_{j_{i}}X_i&+\left(\nabla_{X_i}J\right)\circ d\pi_i\circ j_i\\
&=O\left(\|J_i-J\|_{C^0}.\|X_i\|_{W^{k,p}}\right)+O\left(\|X_i\|_{W^{k,p}}^2\right),
\end{split}
\end{displaymath}
for some fixed values of $k,p$, say $k=p=2$. If $X_i\neq 0$, 
re-scaling $X_i$ we obtain the section $\ov{X_i}$ with 
$\|\ov{X_i}\|_{W^{k,p}}=1$ in the domain of $\Kcal_{\pi_i,J_i}$ such that 
$\|\Kcal_{\pi_i,J_i}(\ov{X_i})\|$ becomes arbitrarily small as $i$ goes to infinity. 
However, this implies that the kernel of 
$\Kcal_{\pi,J}$ is non-trivial, unless $X_i$ are zero for $i$ large enough.
Since $C$ is $(n,K)$-rigid, the later should be the case, i.e. 
the image of $\Sig_i$ is included in $C_i$ for $i$ large enough.
The somewhere injectivity assumption on $f_i:\Sig_i\ra M$ then
implies that $\Sig_i=C_i$.\\
Thus the limit curve $f:\Sig\ra C\subset M$ is of degree $1$ and is an embedding. 
Furthermore, except for finitely many values of $i$, $f_i:\Sig_i\ra M$ is 
an embedding as well. Since $f_i:\Sig_i\ra M$ are critical 
$(J_i,n,K)$ curves, there are branched covering maps $p_i:\Sig_i'\ra \Sig_i$
of bounded degree (bounded by $n$) 
such that $\Ker(\Kcal_{p_i})\neq 0$. Furthermore, the genus of $\Sig_i'$ is 
at most $h+K$, where $h$ denotes the genus of $C$.  
By passing to a further sub-sequence we may assume that 
$\deg(p_i)=d\leq n$, and that the genus of all of $\Sig_i$ is $g\leq K+h$.
We may also assume that the maps $p_i$ converge to a degree $d$ covering map 
$p:\Sig'\ra \Sig$. Since the kernel of each $\Kcal_{p_i}$ is non-trivial,
the kernel of $\Kcal_p$ is non-trivial as well, violating the assumption that 
$J\in\Jcal_n(K)$. This contradiction completes the proof. 
\end{proof}

As long as proposition~\ref{prop:intersection-is-transverse}
is true for a positive integer $n$, the open subsets $\Jcal_n(K)$ are also dense in 
$\Jcal^\infty(M,\omega)$.

\begin{lem}\label{lem:dense}
 The subspace $\Jcal_n(K)\subset \Jcal^\infty(M,\omega)$ is dense with 
respect to $C^\ell$ topology for any $\ell>1$ and $n\leq 4$.
\end{lem}
\begin{proof}
 Suppose that $\ell>1$ and $n\leq 4$. By elliptic regularity
$$\Jcal_n(K)=\Jcal^\infty(M,\omega)\cap \Jcal^\ell_n(K).$$
The argument used in the proof of lemma~\ref{lem:open} 
in the smooth case may be used to 
show that $\Jcal^\ell_n(K)$ is an open subset of 
$\Jcal^\ell(M,\omega)$ with respect to
the $C^\ell$ topology. Moreover, $\Jcal^\ell_n(K)$ contains $\Jcal^\ell_n(M,\omega)$
which is dense in $\Jcal^\ell(M,\omega)$ by theorem~\ref{thm:main}, 
and is thus an open dense subset.
Since $\Jcal_{n}(K)$ is the intersection of $\Jcal^\infty(M,\omega)$ 
with an open dense subset of $\Jcal^\ell(M,\omega)$, the claim of the lemma is 
implied. 
\end{proof}

The above two lemmas imply the following stronger 
version of theorem~\ref{thm:main-introduction}.
\begin{thm}\label{thm:main-smooth}
The subset $\Jcal_n^\infty(M,\omega)\subset \Jcal^\infty(M,\omega)$ 
is of second category for 
$n=1,2,3,4$. 
\end{thm}
\begin{proof}
 Since lemma~\ref{lem:dense} is true for all $\ell$, $\Jcal_n(K)$ is dense in 
$\Jcal^\infty(M,\omega)$ in $C^\infty$ topology and 
$$\Jcal_n(M,\omega)=\bigcap_{K\in \Z^{>0}}\Jcal_n(K)$$
is the intersection of a countable collection of open dense subsets 
of $\Jcal^\infty(M,\omega)$.
\end{proof}

\section{Conclusion and final remarks}\label{sec:conclusion}
Let us fix a generic almost complex structure 
$$J\in\Jcal^\infty_{rigid}(M,\omega):=\Jcal^\infty_4(M,\omega),$$
and let $C\subset M$ be any (smooth) $J$-holomorphic curve of genus 
$h\geq 0$ which is 
$4$-rigid. In this case, $C$ determines an open component of the moduli space
$\Modbar_g(\beta,J)$ for any homology class 
$\beta=d[C]\in\Ht_2(M,\Z)$ with $0<d<5$.
This open component may be identified with the moduli space 
$\Modbar_g(C,d[C])$ of branched 
covering maps of $C$ of total genus $g$ and total degree $d$.
Theorem~\ref{thm:main-smooth} tells us that for a generic $J$ 
as above, for any homology class 
$\alpha\in \Ht_2(M,\Z)$ in $M$ with
$c_1(M).\alpha=0$ and $|\alpha|<5$,  all $J$-holomorphic 
curves representing $\alpha$ are embedded,
$4$-rigid, and they are all disjoint from each-other. 
Thus for any embedded $J$-holomorphic curve $C$ 
in $M$ of genus $h$ and representing
the homology class $\alpha$ as above, and for any $g\geq h$ 
and $0<d<5$ the compactified
moduli space $\Modbar_g(C,d[C])$ 
is a component of the Delign-Mumford compactification 
$\Modbar_g(M,d\alpha;J)$ of $\Mod_g(M,d\alpha;J)$.
The contribution of this component of the
moduli space to the Gromov-Witten invariants comes
from integrating the Euler class of an obstruction bundle over the virtual moduli cycle
associated with the moduli space $\Modbar_g(C,d[C])$.  
The fiber of the obstruction  bundle
$\Upsilon=\Upsilon_g(C,d;J)$ at a point $\pi:\Sig\ra C$ in $\Mod_g(C,d[C])$
is the cokernel of the operator $\Kcal_\pi$. 
After a suitable stabilization of this bundle, we may extend it
as a virtual obstruction bundle over $\Modbar_g(C,d[C])$ (see \cite{LT}, or \cite{FO}).
Denote the Euler class of this obstruction bundle by $\chi_g(C,d;J)$. 
Since the rank of the obstruction bundle and the dimension of the virtual 
moduli cycle are both $4((g-1)-d(h-1))$, $\chi_g(C,d;J)$ may be
integrated against the virtual class 
$\big[\Modbar_g(C,d[C])\big]^{vir}$ to give the rational numbers
\begin{displaymath}
C_g(C,d;J):= \int_{\big[\Modbar_g(C,d[C])\big]^{vir}}\chi_g(C,d;J)\in \Qbb.
\end{displaymath}
The number $C_g(C,d;J)$ is the local contribution of the 
curve $C$ to the Gromov-Witten invariant
$N_g(M,d\alpha)$, and depends not only on the normal 
bundle of the $J$-holomorphic curve $C$, but also on
$\nabla J$ in normal directions, or equivalently, on the 
Nijenhuis tensor $N_J$ associated with $J$ 
along $C$. In the particular case of degree one curves, $C_h(C,1;J)$ is always
equal to $\pm1$, depending on the sign of the spectral flow from the
Dolbeault $\dbar$-operator to $\Kcal_{Id}$. Moreover, 
the argument of \cite{P} may be used
to compute all $C_g(C,1;J)$ for $C\in \Mod_h(\beta;J)$ and $g\geq 0$ 
via the following generating function formula
\begin{displaymath}
\sum_{g\geq h}C_g(C,1;J)\lambda^{2g-2}=(2\sin(\frac{\lambda}{2}))^{2h-2}C_h(C,1;J).
\end{displaymath}
Let us now assume that $C=\CP1$ is a rational curve. Note that the moduli space 
$\Modbar_0(\CP1,d[\CP1])$, after moding out by the 
automorphism group of the domain,
may be identified with $\CP{2d-2}$. In fact, for every map $\pi:\CP1\ra\CP1$ 
of degree $d$, the pre-image of three generic points 
(which may be labeled $0,1$ and $\infty$)
consist of $3d$ points, $d$ in the pre-image of each one of them. 
There are thus $d^3$ triples $(x_0,x_1,x_\infty)$
on $\CP1$ which are mapped to $(0,1,\infty)$ by $\pi$. 
Each triple as above determines a re-parametrization
of the domain of $\pi$ which takes $(0,1,\infty)$ to itself. 
Thus, associated with any degree $d$ map 
in $\CP{2d-2}$, there are $d^3$  elements in $\Modbar_0(C,d[C])$ 
which differ only by automorphisms of the 
domain. The fibers of the obstruction bundle over these 
points are naturally mapped to each other 
under the automorphism. Using the frame-work developed 
by Li and Tian \cite{LT}, we may thus write 
\begin{equation}\label{eq:genus-zero-contribution}
\int_{\left[\Modbar_0(C,d[C])\right]^{vir}} \chi_0(C,d[C];J)=\frac{1}
{d^3}\int_{\PP^{2d-2}}\chi_0'(C,d;J)=:\frac{c(C,d;J)}{d^3},
\end{equation}
where $\chi_0'(C,d;J)\in\Ht^{4d-4}(\CP{2d-2},\Z)$, and $c(C,d;J)$ 
is thus an integer (see \cite{Voicin}).
For $0\neq \alpha\in\Ht_2(M,\Z)$ satisfying $c_1(M).\alpha=0$ 
and $|\alpha|\leq 4$ let
$$e_0(\alpha,d;J):=\sum_{C\in\Mod_0(\alpha;J)} c(C,d;J).$$
Let $\Ical_5$ denote the subspace generated by the formal powers $q^\beta$ over 
$\Q$, where $\beta\in\Ht_2(M,\Z)$
satisfies $c_1(M).\beta=0$ and $|\beta|\geq 5$.
Replacing the computation of equation~\ref{eq:genus-zero-contribution} in the 
generating function for genus 
zero Gromov-Witten invariants we obtain
\begin{equation}\label{eq:genus-zero-GW-invariants}
\begin{split}
 \sum_{\substack{0\neq\beta\in \Ht_2(M,\Z)\\ c_1(M).\beta=0}}N_0(\beta)q^\beta=&
\sum_{\substack{0\neq \alpha \in \Ht_2(M,\Z)\\ c_1(M).\alpha=0}}\ \sum_{d>0}\  
\frac{e_0(\alpha,d;J)}{d^3} q^{d\alpha}
\ \ (\mathrm{modulo}\ \Ical_5).\\
\end{split}
\end{equation}
This presentation of genus zero Gromov-Witten invariants of $(M,\omega)$ is a 
re-statement of conjecture~7.4.5 from \cite{CK} 
for homology classes with divisibility less than $5$.\\

In the integrable case, if $M$ is equipped with a complex structure $J$ and 
$C\subset M$ is an isolated $d$-rigid 
smooth $J$-holomorphic curve of genus $h$, the local contributions 
$C_g(C,d,J)$ are proved by Bryan and Pandharipande
\cite{R-Jim} to be independent of the normal bundle, and the 
complex structure on $C$. They show that these contributions only 
depend on the genus $h$ of $C$ and the integers $g\geq h$ and $d>0$.
We may thus denote them by $C_g(h,d)$. However, when the almost complex 
structure $J\in\Jcal^\infty_{rigid}(M,\omega)$ is not integrable, 
the independence of the local contributions $C_g(C,d;J)$ from $J$ is not clear.
If $\{J_t\}_{t\in[0,1]}$ is a generic path of almost complex structures with 
$J_0,J_1\in\Jcal^\infty_{rigid}(M,\omega)$, we can not guarantee that $\{J_t\}_t$
determines a compact one-dimensional oriented cobordism from 
$\Mod_h(\beta;J_0)$ to $\Mod_h(\beta;J_1)$. The smooth moduli space 
corresponding to such path may have ends 
different from  $\Mod_h(\beta;J_0)\cup\Mod_h(\beta;J_1)$.
One may show, however, that every other limit point 
is in correspondence with an embedded $J_t$-holomorphic
curve $C$, for $t$ in a finite subset of $[0,1]$, and that 
for any such curve $C$, there is a 
degree $2$ branched covering map $\pi:\Sig\ra C$ so that 
the kernel of $\Kcal_\pi$ is  one dimensional. Moreover, if $X$
 is the generator of this kernel and $\sig$ is the natural 
involution of $\Sig$ then $\sig^*X=-X$.
As we vary $J$ through the generic path $\{\Jcal_t\}_t$ 
and pass through one of the above finitely 
many values of $t$, the signed count of points in the moduli space 
$\Mod_h(\beta;J)$ {\emph{may}} change by $1$. Formulating a precise 
wall-crossing formula is not straight-forward, although it is possible in principle 
to do so, following the method used by Taubes \cite{Taubes}
 in complex dimension $2$.\\ 



\begin{thebibliography}{Dillo 83}
\bibitem{R-Jim} Bryan, J., Pandharipande, R.,
BPS-states of curves in Calabi-Yau $3$-folds, \emph{ Geom. Topol.}
5 (2001), 287-318 (electronic).
\bibitem{CK} Cox, D., Katz, S., 
\emph{Mirror Symmetry and Algebraic Geometry}, Math. Surveys and Mono. ,
vol 68, AMS, Providence RI, 1999.
\bibitem{FO} Fukaya, K., Ono, K., Arnold conjecture and Gromov-Witten invariants, 
\emph{Topology } 38 (1999), 933-1048.
\bibitem{Go-Va-1} Gopakumar, R., Vafa, C.,
M theory and topological strings II, (1998), arxiv:hep-th/9812127.
\bibitem{IP1} Ionel, E. N., Parker, T.,
The Gromov invariants of Ruan-Tian and Taubes, \emph{ Math. Res.
Lett.}  4  (1997),  no. 4, 521-532.
\bibitem{IP2} Ionel, E. N., Parker, T.,
Relative  Gromov invariants, \emph{Ann. of Math.}  157  (2003),  45-96.
\bibitem{LT} Li, J., Tian, G., Virtual modli cycles and Gromov-Witten invariants
of general symplectic manifolds, \emph{J. of the Amer. Math. Society},
11, no. 1, (1998), pp. 119-174.
\bibitem{Mc-Sa} McDuff, D., Salamon, D.,
\emph{J-holomorphic curves and quantum cohomology},
University Lecture Series, 6. American Mathematical Society, Providence, RI, 1994.
\bibitem{Oh} Oh, Y., Super-rigidity and finiteness of embedded 
$J$-holomorphic curves
on Calabi-Yau threefolds, \emph{preprint, available at }arXiv:0807.3152v1.
\bibitem{OZ} Oh, Y., Zhu, K., Embedding property of $J$-holomorphic curves 
in Calabi-Yau manifolds for generic $J$, {\emph{Asian J. Math.}}13 (2009), 323-340,
\emph{also available at }arXiv:0805.3581v2.
\bibitem{P}  Pandharipande, R. , Hodge integrals and degenerate contributions,
\emph{Comm. Math. Phys.}, 208 (1999), no. 2, 489-506.
\bibitem{PZ} Pandharipande, R. ,  Zinger, A.,
Enumerative Geometry of Calabi-Yau 5-Folds, 
{\emph{ASPM}} 59 (2010), 239-288,
\emph{also available at}
 arXiv:0802.1640.
\bibitem{R-T2} Ruan, Y., Tian, G.,
Higher genus symplectic invariants and sigma model coupled with
gravity, \emph{ Turkish J. Math.} 20 (1996) no.1, 75-83.
\bibitem{Taubes} Taubes, C.H., Counting pseudo-holomorphic submanifolds in
dimension $4$, \emph{J. Differential Geom.} 44 (1996) no.4, 818-893.
\bibitem{Voicin} Voisin, C., A mathematical proof of a formula of Aspinwall 
and Morrison,
 \emph{Compositio Math.} 104 (1996), no. 2, 135--151.
\end{thebibliography}
\end{document}